\documentclass[11pt]{article}
\usepackage{amsmath}
\usepackage{amsthm}
\usepackage{amsfonts}
\usepackage[usenames,dvipsnames]{xcolor}
\usepackage{enumitem}
\usepackage[T1]{fontenc}
\usepackage[utf8]{inputenc}
\newtheorem{theo}{Theorem}[section]
\newtheorem{lem}[theo]{Lemma}

\newtheorem{rem}[theo]{Remark}
\newtheorem{defi}[theo]{Definition}

\newcommand{\mysection}[1]{\section{#1} \setcounter{equation}{0}}
\newcommand{\proofc}{{\sc Proof} \ }
\newcommand{\be}{\begin{equation} \label}
\newcommand{\ee}{\end{equation}}
\newcommand{\bea}{\begin{eqnarray}\label}
\newcommand{\eea}{\end{eqnarray}}
\newcommand{\bas}{\begin{eqnarray*}}
\newcommand{\eas}{\end{eqnarray*}}
\newcommand{\bit}{\begin{itemize}}
\newcommand{\eit}{\end{itemize}}
\renewcommand{\qed}{\hfill$\Box$ \vskip.2cm}
\newcommand{\nn}{\nonumber}
\newcommand{\R}{\mathbb{R}}
\newcommand{\N}{\mathbb{N}}
\newcommand{\pO}{\partial\Omega}
\newcommand{\bom}{\overline{\Omega}}
\newcommand{\eps}{\varepsilon}
\newcommand{\dist}{{\rm dist} \, }
\newcommand{\supp}{{\rm supp} \, }

\newcommand{\wto}{\rightharpoonup}

\newcommand{\io}{\int_\Omega}

\newcommand{\abs}{\\[5pt]}
\newcommand{\ueps}{u_\eps}
\newcommand{\veps}{v_\eps}

\newcommand{\qm}{q_-(p)}
\newcommand{\qp}{q_+(p)}

\usepackage{newunicodechar}
\newunicodechar{ε}{\varepsilon}
\newunicodechar{χ}{\chi}
\newunicodechar{φ}{\varphi}
\newunicodechar{∞}{\infty}
\newunicodechar{ℝ}{\mathbb{R}}
\newunicodechar{ψ}{\psi}
\newunicodechar{ζ}{\zeta}
\newunicodechar{Δ}{\Delta}
\newunicodechar{τ}{\tau}
\newunicodechar{η}{\eta}
\newunicodechar{δ}{\delta}
\newunicodechar{λ}{\lambda}
\newunicodechar{γ}{\gamma}
\newunicodechar{ℕ}{\mathbb{N}}
\newunicodechar{ν}{\nu}
\newunicodechar{∂}{\partial}
\newcommand{\f}[2]{\frac{#1}{#2}}
\newcommand{\na}{\nabla}
\newcommand{\ue}{u_{ε}}
\newcommand{\ve}{v_{ε}}

\newcommand{\Om}{\Omega}
\newcommand{\Ombar}{\overline{\Omega}}
\newcommand{\Lom}[1]{L^{#1}(\Om)}
\newcommand{\LomT}[1]{L^{#1}(\Om\times(0,T))}
\newcommand{\norm}[2][]{\|#2\|_{#1}}
\newcommand{\ptilde}{\widetilde{p}}
\newcommand{\qtilde}{\widetilde{q}}

\newcommand{\intnT}{\int_0^T}
\newcommand{\intnt}{\int_0^t}
\newcommand{\intninf}{\int_0^\infty}
\newcommand{\dOm}{\partial\Omega}
\newcommand{\amrand}{\big|_{\dOm}}
\newcommand{\amrandalltime}{\big|_{\dOm\times(0,\infty)}}
\newcommand{\intdom}{\int_{\dOm}}
\newcommand{\set}[1]{\{#1\}}
\DeclareMathOperator{\essinf}{ess inf}
\DeclareMathOperator{\arcoth}{arcoth}
\DeclareMathOperator{\diam}{diam}

\newcommand{\kl}[1]{\left(#1\right)}
\newcommand{\kkl}[1]{\left\{#1\right\}}
\newcommand{\sub}{\subset}
\allowdisplaybreaks
\begin{document}
\title{A generalized solution concept for the Keller-Segel system with logarithmic sensitivity: 
Global solvability for large nonradial data}
%
%
%
\author{
Johannes Lankeit\footnote{jlankeit@math.uni-paderborn.de}\\
{\small Institut f\"ur Mathematik, Universit\"at Paderborn,}\\
{\small 33098 Paderborn, Germany} 
\and
Michael Winkler\footnote{michael.winkler@math.uni-paderborn.de}\\
{\small Institut f\"ur Mathematik, Universit\"at Paderborn,}\\
{\small 33098 Paderborn, Germany} }
%
%
%
\maketitle
\begin{abstract}
\noindent 
  The chemotaxis system 
  \bas
	\left\{ \begin{array}{l}
	u_t = \Delta u - \chi\nabla  \cdot (\frac{u}{v}\nabla v), \\[1mm]
	v_t=\Delta v - v+u,
	\end{array} \right.
  \eas
  is considered in a bounded domain $\Omega\subset \R^n$ with smooth boundary, where $\chi>0$.\abs
  An apparently novel type of generalized solution framework is introduced within which 
  an extension of previously known ranges for the key parameter $\chi$ with regard to global solvability is achieved.
  In particular, it is shown that under the hypothesis that 
  \bas
	\chi<\left\{ \begin{array}{ll}
	\infty
	\qquad & \mbox{if } n=2, \\[1mm]
	\sqrt{8}
	\qquad & \mbox{if } n=3, \\[1mm]
	\frac{n}{n-2}
	\qquad & \mbox{if } n\ge 4,
	\end{array} \right.
  \eas
  for all initial data satisfying suitable assumptions on regularity and positivity,
  an associated no-flux initial-boundary value problem admits a globally defined generalized solution.
  This solution inter alia has the property that
  \bas
	u\in L^1_{loc}(\bom\times [0,\infty)).
  \eas
  {\bf Math Subject Classification (2010): 35K55, 35D99, 92C17} 
\end{abstract}
%
%
%


%
%
\newpage
\section{Introduction}\label{intro}
We consider the Keller-Segel system with logarithmic sensitivity, as given by the initial-boundary value problem
\be{01}
	\left\{ \begin{array}{ll}
	u_t=\Delta u -\chi \nabla\cdot \Big(\frac{u}{v}\nabla v\Big), 
	\qquad & x\in\Omega, \ t>0, \\[1mm]
	\frac{\partial u}{\partial\nu}=0,
	\qquad & x\in\pO, \ t>0, \\[1mm]
	u(x,0)=u_0(x),
	\qquad x\in\Omega,
	\end{array} \right.
\ee
coupled to the parabolic problem
\be{02}
	\left\{ \begin{array}{ll}
	v_t=\Delta v-v+u,
	\qquad & x\in\Omega, \ t>0, \\[1mm]
	\frac{\partial v}{\partial\nu}=0,
	\qquad & x\in\pO, \ t>0, \\[1mm]
	v(x,0)=v_0(x),
	\qquad x\in\Omega,
	\end{array} \right.
\ee
where $\Omega$ is a bounded domain in $\R^n$, $n\ge 2$, with smooth boundary, $\chi$ is a positive parameter
and the given initial data $u_0$ and $v_0$ satisfy suitable regularity and positivity assumptions.\abs
This system can be viewed as a prototypical parabolic model for self-enhanced chemotactic migration processes
in which cross-diffusion 
occurs in accordance with the Weber-Fechner law of stimulus perception (\cite{hillen_painter2009},
\cite{rosen}), and accordingly a considerable literature is concerned with its mathematical analysis.
However, up to now it seems yet unclear to which extent the particular mechanism of
taxis inhibition at large signal densities in (\ref{01}) is sufficient to prevent phenomena of blow-up,
known as the probably most striking qualitative feature of the classical Keller-Segel system: 
Indeed, in its fully parabolic version, as determined by the choice $\tau:=1$ in
\be{KS}
	\left\{ \begin{array}{l}
	u_t=\Delta u - \chi \nabla \cdot (u\nabla v), \\[1mm]
	\tau v_t=\Delta v-v+u,
	\end{array} \right.
\ee
the latter admits solutions blowing up in finite time for any choice of $\chi>0$ whenever $n\ge 2$
(\cite{herrero_velazquez}, \cite{win_JMPA}), and in the simplified parabolic-elliptic case obtained on choosing
$\tau:=0$ it is even known that some radial solutions to an associated Cauchy problem in the whole plane collapse into
a persistent Dirac-type singularity in the sense that a globally defined measure-valued solution exists which 
has a singular part beyond some finite time and asymptotically approaches a Dirac measure 
(cf.~e.g.~\cite{tello_win_Pisa} or also \cite{luckhaus_sugiyama_ARMA}).\abs
As opposed to this, the literature has identified various circumstances under which phenomena of this type are ruled
out in (\ref{01})-(\ref{02}): For instance, when $\chi<\chi_0(n)$ with some $\chi_0(2)>1.015$ and 
$\chi_0(n):=\sqrt{\frac{2}{n}}$ for $n\ge 3$,
global bounded classical solutions exist
for all reasonably regular positive initial data (\cite{lankeit}, \cite{biler}, \cite{fujie_bounded}, 
\cite{zhao_zheng_sining_JMAA2016}, \cite{mizukami_yokota}, \cite{win_MMAS});
in the corresponding parabolic-elliptic analogue, the same conclusion holds with $\chi_0(2)=\infty$
(\cite{fujie_senba}) and with $\chi_0(n):=\frac{2}{n-2}$ when $n\ge 3$ and the spatial setting is radially
symmetric (\cite{nagai_senba1998}, cf.~also \cite{fujie_senba_NON} for a related result addressing
a variant with its second equation being $\tau v_t=\Delta v - v+u$ for small $\tau>0$),
whereas it is known that some exploding solutions exist if $n\ge 3$
and $\chi>\frac{2n}{n-2}$ (\cite{nagai_senba1998}).
As for larger values of $\chi$ in the fully parabolic problem (\ref{01})-(\ref{02}), 
in some cases		
at least certain global generalized solutions can be found which satisfy 
\be{l1}
	u\in L^1_{loc}(\bom\times [0,\infty))
\ee
and thereby indicate the absence of strong singularity formation of the flavor described above.
Such constructions are possible in the context of natural weak solution concepts if 
\be{cond1}
	\chi<\sqrt{\frac{n+2}{3n-4}}
\ee
(\cite{win_MMAS}) and within a slightly more generalized framework if merely 
\be{cond2}
	\chi<\sqrt{\frac{n}{n-2}}
\ee
but in addition
the solutions are supposed to be radially symmetric (\cite{stiwi}).
To the best of our knowledge, however, the question how far (\ref{cond1}) is optimal with respect to the existence 
of not necessarily radial solutions fulfilling (\ref{l1}) is yet unsolved;
in particular, it appears to be unknown whether in nonradial planar settings such solutions do exist also
beyond the range $\chi<\sqrt{2}$ determined by (\ref{cond1}).\abs
{\bf Main results.} \quad
The purpose of this work is to design a novel concept of generalized solvability which is 
yet suitably strong so as to require (\ref{l1}), but which on the other hand is mild enough so that it 
enables us to construct corresponding global solutions without any symmetry hypotheses and under conditions
somewhat weaker than (\ref{cond1}) and actually also than (\ref{cond2}). 
More precisely, considering (\ref{01})-(\ref{02}) under the assumptions that
\be{init}
    \left\{
    \begin{array}{ll}
    u_0 \in C^0(\bom) \; &\mbox{ is such that } u_0\ge 0 \mbox{ in $\Omega$ and }
	u_0 \not\equiv 0, \quad
	\mbox{and that} \\
    v_0 \in W^{1,\infty}(\Omega) \; &\mbox{ satisfies $v_0> 0$ in $\bom$,}
    \end{array}
    \right.
\ee
we can state our main results as follows.
\begin{theo}\label{theo14}
  Let $n\ge 2$ and $\Omega\subset\R^n$ be a bounded domain with smooth boundary, and let $\chi>0$ be such that
  \be{14.1}
	\chi<\left\{ \begin{array}{ll}
	\infty
	\qquad & \mbox{if } n=2, \\[1mm]
	\sqrt{8}
	\qquad & \mbox{if } n=3, \\[1mm]
	\frac{n}{n-2}
	\qquad & \mbox{if } n\ge 4.
	\end{array} \right.
  \ee
  Then for any $u_0$ and $v_0$ fulfilling (\ref{init}), the problem (\ref{01})-(\ref{02}) possesses at least
  one global generalized solution $(u,v)$ in the sense of Definition \ref{defi12} below. 
  In particular, this solution satisfies (\ref{l1}), and moreover we have
  \be{mass_u}
	\io u(\cdot,t)=\io u_0
	\qquad \mbox{for a.e.~} t>0.
  \ee
\end{theo}
{\bf Plan of the paper.} \quad
Our approach will be based on the essentially well-known fact that the functional $\io u^p v^q$
enjoys certain quasi-entropy features along trajectories of (\ref{01})-(\ref{02}), provided that the crucial positive
parameter $p$ therein satisfies $p<\frac{1}{\chi^2}$ (cf.~Section \ref{sect4} for a corresponding observation 
addressing global smooth solutions to the regularized problems (\ref{0eps}) below). 
The challenge now consists in taking appropriate advantage of correspondingly implied a priori estimates
obtained in Sections \ref{sect5}, \ref{sect6} and \ref{sect7}, which inter alia
seem far from sufficient to warrant $L^1$ bounds for the cross-diffusive flux $\chi \frac{u}{v}\nabla v$ 
especially in cases when $\chi$ is large and hence $p$ needs to be chosen small.\\ 
In the preparatory Section \ref{sect3}, 
we will therefore resort to a solution framework involving certain sublinear powers of $u$ rather than $u$ itself,
thus reminiscent of the celebrated concept of renormalized solutions {\cite{diperna_lions}). 
This idea has partially been adapted to the present context in \cite{stiwi} already, but in the present work we shall
further weaken the requirements on solutions to a considerable extent: Namely, for the crucial first sub-problem
(\ref{01}) to be solved we shall only require that the {\em coupled} quantity $u^p v^q$, with certain positive
$p$ and $q$, satisfies a 
parabolic {\em inequality} associated with (\ref{01})-(\ref{02}) in a weak form, 
and that moreover $\io u(\cdot,t) \le \io u_0$ for a.e.~$t>0$; a key observation, to be made in Lemma \ref{lem13}, 
will reveal that if we furthermore assume the component $v$ to fulfill (\ref{02}) in a natural weak sense, then
we indeed obtain a concept consistent with that of classical solvability in (\ref{01})-(\ref{02})
for all suitably smooth functions.\\
As seen in Section \ref{sect8} by means of appropriate compactness arguments, the previously gained estimates
in fact enable us to construct a global solution within this framework.
\mysection{A concept of generalized solvability}\label{sect2}
In specifying the subsequently pursued concept of weak solvability, we 
first require certain products $u^p v^q$ to satisfy an inequality
which can be viewed as generalizing a classical supersolution property of this quantity with regard to 
(\ref{01})-(\ref{02}). 
\begin{defi}\label{defi9}
  Let $p\in (0,1)$ and $q\in (0,1)$, and suppose that $u$ and $v$ are measurable functions on $\Omega\times (0,\infty)$
  such that $u> 0$ and $v>0$ a.e.~in $\Omega\times (0,\infty)$, that 
  \be{9.1}
	u^p v^q \in L^1_{loc}(\bom\times [0,\infty))
	\qquad \mbox{and} \qquad
	u^{p+1} v^{q-1} \in L^1_{loc}(\bom\times [0,\infty)),
  \ee
  and that $\nabla u^\frac{p}{2}$ and $\nabla v^\frac{q}{2}$ belong to $L^1_{loc}(\Omega\times (0,\infty))$ 
  and are such that
  \be{9.2}
	v^\frac{q}{2}\nabla u^\frac{p}{2} \in L^2_{loc}(\bom\times [0,\infty))
	\qquad \mbox{and} \qquad
	u^\frac{p}{2}\nabla v^\frac{q}{2} \in L^2_{loc}(\bom\times [0,\infty)).
  \ee
  Then $(u,v)$ will be called a {\em global weak $(p,q)$-supersolution} of (\ref{01}) if 
  \bea{9.3}
	& & \hspace*{-20mm}
	- \int_0^\infty \io u^p v^q \varphi_t
	- \io u_0^p v_0^q \varphi(\cdot,0) \nn\\
	&\ge& \frac{4(1-p) q - 4q^2 - p(1-p)^2 \chi^2}{pq(p\chi+1-q)} 
		\int_0^\infty \io v^q |\nabla u^\frac{p}{2}|^2 \varphi \nn\\
	& & + \frac{4(pχ+1-q)}q
		\int_0^\infty \io \Big| u^{\f p2}\na v^{\f q2} 
		- \frac{(1-p)χ+2q}{2(pχ+1-q)} v^{\f q2}\na u^{\f p2}\Big|^2 \varphi \nn\\
	& & - \frac{2p\chi}{q} \int_0^\infty \io u^\frac{p}{2} v^q \nabla u^\frac{p}{2} \cdot\nabla\varphi \nn\\
	& & + \Big(1-\frac{p\chi}{q}\Big) \int_0^\infty \io u^p v^q \Delta\varphi \nn\\
	& & - q\int_0^\infty \io u^p v^q \varphi
	+ q\int_0^\infty \io u^{p+1} v^{q-1} \varphi
  \eea
  for all nonnegative $\varphi\in C_0^\infty(\bom\times [0,\infty))$ such that $\frac{\partial\varphi}{\partial\nu}=0$
  on $\pO\times (0,\infty)$ and if moreover 
 \begin{align}
  u^pv^q &> 0 \qquad \text{a.e. on } \dOm\times(0,\infty)\label{eq:positivityonbdry}.
 \end{align}
\end{defi}
\begin{rem}\label{rem:lsgdef}
\quad (i) 
  Observing that (\ref{9.1}) in particular ensures that 
  $u^\frac{p}{2} v^\frac{q}{2} \in L^2_{loc}(\bom\times [0,\infty))$, and that hence (\ref{9.1}) and (\ref{9.2}) 
  warrant that
  \bas
	u^\frac{p}{2} v^q \nabla u^\frac{p}{2} 
	= (u^\frac{p}{2} v^\frac{q}{2}) \cdot (v^\frac{q}{2} \nabla u^\frac{p}{2})
	\in L^1_{loc}(\bom\times [0,\infty))
  \eas
  and similarly $u^p v^\frac{q}{2}\nabla v^\frac{q}{2}\in L^1_{loc}(\bom\times [0,\infty))$, it follows that
  under the above requirements all integrals in (\ref{9.3}) are indeed well-defined.\\
 (ii) According to \eqref{9.1} and \eqref{9.2}, for a.e. $t>0$,  $u^{\f p2}(\cdot,t)v^{\f q2}(\cdot,t)\in W^{1,2}(\Om)$ so that $u^{\f p2}v^{\f q2}(\cdot,t)\big|_{\dOm}\in L^2(\dOm)$ exists in the sense of traces, giving meaning to the positivity requirement in \eqref{eq:positivityonbdry}.
\end{rem}
Apart from that, we will require the second problem (\ref{02}) to be satisfied in the following rather natural 
weak sense.
\begin{defi}\label{defi11}
  A pair $(u,v)$ of functions
  \be{11.1}
	\left\{ \begin{array}{l}
	u\in L^1_{loc}(\bom\times [0,\infty)), \\[1mm]
	v\in L^1_{loc}([0,\infty);W^{1,1}(\Omega))
	\end{array} \right.
  \ee
  will be named a {\em global weak solution} of (\ref{02}) if
  \be{11.2}
	- \int_0^\infty \io v\varphi_t - \io v_0 \varphi(\cdot,0)
	= - \int_0^\infty \io \nabla v\cdot\nabla\varphi
	- \int_0^\infty \io v\varphi
	+ \int_0^\infty \io u\varphi
  \ee
  for all $\varphi\in C_0^\infty(\bom\times [0,\infty))$.
\end{defi}
Following an approach already pursued in \cite{win_SIMA} in a considerably less involved related context, 
in order to complete our solution concept we will complement the above two requirements by merely 
postulating an upper bound for the mass functional $\io u$ in terms of $\io u_0$:
\begin{defi}\label{defi12}
  A couple of nonnegative measurable functions $u$ and $v$ defined on $\Omega\times (0,\infty)$ will be said to be
  a {\em global generalized solution} of (\ref{01})-(\ref{02}) if $(u,v)$ is a global weak solution of (\ref{02})
  according to Definition \ref{defi11}, if there exist $p\in (0,1)$ and $q\in (0,1)$ such that $(u,v)$ is
  a global weak $(p,q)$-supersolution of (\ref{01}) in the sense of Definition \ref{defi9}, and if moreover
  \be{12.1}
	\io u(\cdot,t) \le \io u_0
	\qquad \mbox{for a.e.~} t>0.
  \ee
\end{defi}
%
%
%
%
%
%
%
This is indeed consistent with the concept of classical solvability in the following sense.
\begin{lem}\label{lem13}
  Let $\chi>0$, and suppse that $(u,v) \in (C^0(\bom\times [0,\infty))\cap C^{2,1}(\bom\times (0,\infty)))^2$
  is such that $(u,v)$ is a global generalized solution
  of (\ref{01})-(\ref{02}) in the sense of Definition \ref{defi12}. Then $(u,v)$ satisfies (\ref{01})-(\ref{02})
  classically in $\Omega\times (0,\infty)$.
\end{lem}
\proof
By means of standard arguments relying on the assumed regularity properties of $v$, it can easily
be verified that $v$ solves \eqref{02} classically.
According to the maximum principle, $v$ hence is strictly positive in $\Ombar\times[0,∞)$ and $v^{q-1}$ is uniformly bounded in every set $\Om\times[0,T)$ for $T\in (0,∞)$. Positivity of $v$ ensures that by \eqref{eq:positivityonbdry} $u>0$ on a dense 
subset of $\dOm\times(0,\infty)$ which moreover is open in $\pO\times (0,\infty)$ by continuity of $u$.

For arbitrary $ψ\in C^{∞}(\Ombar)$ with $ψ\geq 0$ and $\f{∂ψ}{∂ν}\amrand=0$, testing \eqref{9.3} by $φ(x,t):=ψ(x)(1-\frac1{ε}t)_+$, $ε\in(0,1)$, which is permissible by Weierstrass' theorem, and invoking Lebesgue's dominated convergence theorem and continuity of $t\mapsto \io u^p(\cdot,t)v^q(\cdot,t)$ at $t=0$ in taking $ε\searrow 0$ we readily achieve 
\[
 \io u^p(\cdot,0)v^q(\cdot,0)ψ \geq \io u_0^pv_0^q ψ \qquad \text{ for all } ψ\in C^{∞}(\Ombar), ψ\geq 0, \f{∂ψ}{∂ν}\amrand=0,
\]
showing that $u^p(\cdot,0)v^q(\cdot,0)\geq u_0^pv_0^q$ throughout $\Om$. Because of $v(\cdot,0)=v_0>0$ and the monotonicity of $(\cdot)^{\f 1p}$ we obtain $u(\cdot,0)\geq u_0$ in $\Om$ and from continuity of $u$ and \eqref{12.1} we can conclude that $u(\cdot,0)=u_0$ in $\Om$. 

In the first two integrals on the right of \eqref{9.3} straightforward computations yield 
\begin{align}\label{eq:rearrange}
\frac{4(1-p)}{p}& v^q|\na u^{\frac p2}|^2 - \Big(\frac{4(1-p)}{q}+8\Big) u^{\f p2}v^{\f q2} \na u^{\f p2}\na v^{\f q2}+ \frac{4(pχ+1-q)}q u^p|\na v^{\f q2}|^2\nn\\
 &=\frac{4(pχ+1-q)}q \kkl{ u^p|\na v^{\f q2}|^2 - \frac{(1-p)χ+2q}{pχ+1-q} u^{\f p2}v^{\f q2} \na u^{\f p2}\na v^{\f q2}+\frac{(χ-pχ+2q)^2}{4(pχ+1-q)^2}v^q|\na u^{\f p2}|^2}\nn\\
 &+\kkl{\frac{4(1-p)}p-\frac{((1-p)χ+2q)^2}{q(pχ+1-q)}} v^q|\na u^{\f p2}|^2\nn\\
 &= \frac{4(pχ+1-q)}{q} \bigg|u^{\f p2}\na v^{\f q2}- \frac{(1-p)χ+2q}{2(pχ+1-q)}v^{\f q2}\na u^{\f p2}\bigg|^2 +\frac{4(1-p)q-4q^2-p(1-p)^2χ^2}{pq(pχ+1-q)} v^q|\na u^{\f p2}|^2,
\end{align}
since 
\bas
	\frac{4(1-p)}p - \frac{((1-p)χ+2q)^2}{q(pχ+1-q)} 
	&=& \frac{4(1-p)q(pχ+1-q)-((1-p)χ+2q)^2p}{pq(pχ+1-q)} \\
	&=& \frac{4(1-p)q-4q^2-p(1-p)^2χ^2}{pq(pχ+1-q)}.
\eas
In preparation of the following calculations we also note that for each positive function $w\in C^2(\Ombar)$ and any $r>0$, we have the pointwise identities 
\begin{align}\label{eq:deltaup2}
 w^{\f r2}\Delta w^{\f r2} &= w^{\f r2} \na \cdot\Big(\f r2 w^{\f{r-2}2}\na w\Big) 
 =w^{\f r2}\Big(\frac{r(r-2)}4 w^{\f {r-4}2}|\na w|^2+\f r2 w^{\f {r-2}2}\Delta w\Big)\nn\\
&=\frac{r(r-2)}4 w^{r-2}|\na w|^2 +\f r2 w^{r-1}\Delta w =\f{r-2}r |\na w^{\f r2}|^2+\f r2 w^{r-1}Δw
\end{align}
and
\begin{align}\label{eq:deltaup}
 \Delta w^r =\na \cdot (rw^{r-1}\na w) = r(r-1) w^{r-2} |\na w|^2 + rw^{r-1}\Delta w = \f{4(r-1)}r|\na w^{\f r2}|^2 +rw^{r-1}Δw 
\end{align}
The positivity requirement on $w$ in \eqref{eq:deltaup2} and \eqref{eq:deltaup} prompts us to perform the following calculations only for test functions $φ$ compactly supported in $\set{u>0}:=\set{(x,t)\in\Ombar\times[0,∞):\, u(x,t)>0}$, ensuring strict positivity of $u$ and boundedness of $u^{p-1}$ on $\supp φ$.  

Accordingly, for all nonnegative $φ\in C_0^\infty(\Ombar\times(0,∞))$ with $\supp φ\subset \set{u>0}$ and 
$\f{∂φ}{∂ν}\amrandalltime=0$, by \eqref{eq:deltaup} applied to $u$ and $p$, an integration by parts in the integral 
in \eqref{9.3} containing $\na φ$ yields
\begin{align}\label{eq:intbypartsnaphi}
 -\frac{2pχ}{q} &\intninf\io u^{\f p2}v^q \na u^{\f p2}\na φ = \f{2pχ}q \intninf \io v^q|\na u^{\f p2}|^2φ+\f{4pχ}q \intninf \io u^{\f p2}v^{\f q2} \na v^{\f q2}\na u^{\f p2} φ \nn\\
&\qquad + \f{2pχ}q \intninf u^{\f p2}v^q \Delta u^{\f p2}φ - \f{2pχ}q \intninf\intdom u^{\f p2}v^q \f{∂u^{\f p2}}{∂ν} φ\nn\\
&= \f{2pχ}q \intninf \io v^q|\na u^{\f p2}|^2φ+\f{4pχ}q \intninf \io u^{\f p2}v^{\f q2} \na v^{\f q2}\na u^{\f p2} φ \nn\\
&+\f{2(p-2)χ}q \intninf \io v^q|\na u^{\f p2}|^2φ + \f{p^2χ}q\intninf \io u^{p-1}v^q Δu φ- \f{2pχ}q \intninf\intdom u^{\f p2}v^q \f{∂u^{\f p2}}{∂ν} φ,
\end{align}
whereas integrating by parts twice in the integral containing $Δφ$ in \eqref{9.3}, by \eqref{eq:deltaup} applied to $u,p$ and $v,q$, respectively, leads to 
\begin{align}\label{eq:intbypartsdeltaphi}
 \Big(1-\f{pχ}q\Big)\intninf\io u^pv^qΔφ &=  
 \Big(1-\f{pχ}q\Big)\intninf\io v^q\Delta(u^p)φ + \Big(1-\f{pχ}q\Big)\intninf\io u^p\Delta(v^q)φ \nn\\
&+ 2\Big(1-\f{pχ}q\Big)\intninf\io 2u^{\f p2}\na u^{\f p2}2v^{\f q2}\na v^{\f q2}φ\nn\\
&-\Big(1-\f{pχ}q\Big)\intninf\intdom 2u^{\f p2}\f{∂u^{\f p2}}{∂ν}v^qφ-2\Big(1-\f{pχ}q\Big)\intninf\intdom u^p v^{\f q2} \f{∂v^{\f q2}}{∂ν}φ \nn\\
 &\hspace{-2cm}= \f{4(p-1)}p\Big(1-\f{pχ}q\Big) \intninf\io v^q|\na u^{\f p2}|^2φ + \f{4(q-1)}q\Big(1-\f{pχ}q\Big)\intninf\io u^p|\na v^{\f q2}|^2φ\nn \\
  &+\Big(1-\f{pχ}q\Big)p\intninf\io v^q u^{p-1} Δuφ + \Big(1-\f{pχ}q\Big)q \intninf \io u^pv^{q-1}Δvφ\nn\\
 &+ 8\Big(1-\f{pχ}q\Big)\intninf\io u^{\f p2}v^{\f q2} \na u^{\f p2}\na v^{\f q2}φ\nn\\
&-\Big(1-\f{pχ}q\Big)\intninf\intdom 2u^{\f p2}\f{∂u^{\f p2}}{∂ν}v^qφ
\end{align}
for any such $\varphi$,
for we already know that $\f{∂v}{∂ν}=0$ on $\dOm\times(0,∞)$. 
If we combine \eqref{9.3} with \eqref{eq:rearrange}, \eqref{eq:intbypartsnaphi} and \eqref{eq:intbypartsdeltaphi}, we obtain 
\begin{align}\label{eq:ueqtested}
\intninf\io (u^pv^q)_tφ
&\geq 
\kkl{\f{4(1-p)}p + \f{2pχ}q + \f{4(p-1)}p \Big(1-\f{pχ}q\Big) + \f{2(p-2)χ}{q}} \intninf\io v^q|\na u^{\f p2}|^2φ\nn\\
&+\kkl{\f{4(pχ+1-q)}q+\f{4(q-1)}q \Big(1-\f{pχ}q\Big)} \intninf \io u^p|\na v^{\f q2}|^2φ\nn\\
&+\kkl{-\f{4(1-p)χ}q - 8 + \f{4pχ}q +8 - \f{8pχ}q}  \intninf \io u^{\f p2} v^{\f q2} \na u^{\f p2}\na v^{\f q2}φ\nn\\
&+\kkl{\f{p^2χ}q +\Big(1-\f{pχ}q\Big)p} \intninf \io v^q u^{p-1}\Delta u φ \nn\\
&+\Big(1-\f{pχ}q\Big)q\intninf \io u^pv^{q-1}Δvφ - q\intninf\io u^pv^qφ +q\intninf\io u^{p+1}v^{q-1} φ \nn\\
&- \f{2pχ}q \intninf\intdom u^{\f p2}v^q \f{∂u^{\f p2}}{∂ν} φ-\kl{1-\f{pχ}q}\intninf\intdom 2u^{\f p2}\f{∂u^{\f p2}}{∂ν}v^qφ\nn\\
&= 
\intninf\io \kkl{\f{4pχ}{q^2} u^p|\na v^{\f q2}|^2 - \f{4χ}q u^{\f p2}v^{\f q2}\na u^{\f p2}\na v^{\f q2} + pv^qu^{p-1}Δu -pχ u^pv^{q-1}Δv}φ\nn\\
&+q\intninf\io u^pv^{q-1} \kkl{ Δv - v + u }φ \nn\\
&-2\intninf\intdom u^{\f p2}\f{∂u^{\f p2}}{∂ν}v^q φ 
\end{align}
for every $φ\in C_0^\infty(\Ombar\times(0,∞))$ satisfying $φ\geq 0$ throughout $\Om\times(0,\infty)$ and $\f{∂φ}{∂ν}\amrandalltime=0$ as well as $\supp φ\subset\set{u>0}$. 

The observations that 
\begin{align*}
 pu^{p-1}v^q\Big(Δu - χ\na\cdot(\f uv\na v)\Big) &= pu^{p-1}v^qΔu -pχu^{p-1}v^{q-1}\na u\na v + pχu^pv^{q-2}|\na v|^2 - pχ u^{p-1}v^{q-1}Δv\\
 &=pu^{p-1}v^qΔu-\f{4χ}qu^{\f p2}v^{\f q2}\na u^{\f p2}\na v^{\f q2}+ \f{4pχ}{q^2} u^p|\na v^{\f q2}|^2- pχu^{p-1}v^{q-1}Δv, 
\end{align*}
and that $v$ solves \eqref{02}, now turn \eqref{eq:ueqtested} into 
\begin{equation}\label{eq:utestedfinal}
 p \intninf \io u^{p-1}v^q u_t φ \geq p \intninf \io u^{p-1}v^q\kkl{Δu-χ\na\cdot\kl{\f uv\na v}}φ - p \intninf\intdom u^{p-1}v^q \f{∂u}{∂ν} φ
\end{equation}
for all nonnegative $φ\in C_0^\infty(\Ombar\times(0,∞))$ with $\supp φ\subset\set{u>0}$ and $\f{∂φ}{∂ν}\amrandalltime=0$.

Specializing this to nonnegative $φ\in C_0^\infty(\Om\times (0,∞)\cap \set{u>0})$ 
by a Du Bois-Reymond lemma type argument we conclude 
\begin{equation}\label{eq:utineq}
 u_t\geq Δu-χ\na \cdot\kl{\f uv\na v} \qquad \text{in } \set{u>0}.
\end{equation}
Density of $\set{u>0}$ in $\Om\times(0,∞)$, obtained from the assumption that $u>0$ a.e., and continuity show that \eqref{eq:utineq} actually holds on all of $\Om\times(0,\infty)$. 

We pick $t_0>0$ and some nonnegative $ψ\in C^1(\Ombar)$ with $\f{∂ψ}{∂ν}\amrand=0$ such that $\suppψ\subset\set{u(\cdot,t_0)>0}:=\set{x\in \Ombar:\, u(x,t_0)>0}$. Then 
by continuity of $u$ we can find some $τ>0$ such that $\supp ψ\subset \cap_{t\in(t_0-τ,t_0+τ)} \set{u(\cdot,t)>0}$. Applying \eqref{eq:utestedfinal} to functions of the form $φ(x,t)=ζ(t)ψ(x)$, $ζ\in C_0^\infty((t_0-τ,t_0+τ))$ by 
once more invoking a Weierstrass type density argument and the Du Bois-Reymond lemma, we see that 
\[
 \io u^{p-1}v^q u_t(\cdot,t) ψ \geq \io u^{p-1}v^q\kkl{Δu-χ\na\cdot\kl{\f uv\na v}}(\cdot,t) ψ-\intdom u^{p-1}v^q\f{∂u(\cdot,t)}{∂ν} ψ 
\]
for every nonnegative $ψ\in C^1(\Ombar)$ such that $\f{∂ψ}{∂ν}\amrand=0$, $\suppψ\subset\set{u(\cdot,t_0)>0}$ and for almost every $t\in(t_0-τ,t_0+τ)$ -- and due to continuity especially for $t=t_0$. 
In particular inserting 
$ψ_{ε}(x):=(1-\frac 1{ε} \dist(x,\dOm))_+\cdot ψ$ and Lebesgue's theorem show that for every $t>0$, $ψ\in C^1(\Ombar)$, $ψ\geq 0$ with $\f{∂ψ}{∂ν}\amrand=0$ and $\supp ψ\cap \dOm\subset\set{u(\cdot,t)>0}$,  
\begin{equation}\label{eq:utestedbdryfinal}
 \intdom u^{p-1}v^q\f{∂u}{∂ν} ψ\geq 0. 
\end{equation}
Since the integral only depends on $ψ\amrand$ and not on the values of $ψ$ inside $\Omega$, it can be seen that \eqref{eq:utestedbdryfinal} actually holds for any $t>0$ and any nonnegative $ψ\in C^0(\Ombar)$ such that $\supp ψ\cap \dOm\subset\set{u(\cdot,t)>0}$. 

If $\f{∂u}{∂ν}(x_0,t_0)<0$ for some $(x_0,t_0)\in \dOm\times(0,\infty)$ with $u(x_0,t_0)>0$, we pick $\psi_1\in C^0(\dOm)$ such that $\psi_1(x_0)>0$, $\psi_1 \geq 0$, and $\supp \psi_1\subset \{x\in\dOm: \f{∂u}{∂ν}(x,t_0)<0, u(x,t_0)>0 \}=:M$. Moreover, we let $d:=\dist(\supp \psi_1,\dOm\setminus M)$ (or $d=1$ if $\dOm\setminus M=\emptyset$) and let $ψ_2$ be the solution to $-Δψ_2=0$ in $\Om$, $ψ_2=-u^{1-p}(\cdot,t_0)v^{-q}(\cdot,t_0)\psi_1\f{∂u}{∂ν}$ on $\dOm$. Then, thanks to the choice of $\supp ψ_1$, $ψ_2$ is nonnegative on the boundary and hence by the maximum principle in $\Om$. Defining $ψ_3(x):=ψ_2(x)(1-\f2d \dist(x,M))_+$ 
we obtain a nonnegative continuous function on $\Ombar$ whose support intersects the boundary only in $\set{x\in\dOm:\, u(x,t)>0}$ and which hence is a permissible test function in \eqref{eq:utestedbdryfinal}. We conclude that $0\leq \intdom u^{p-1}v^q\f{∂u}{∂ν} ψ_3 =-\intdom |\f{∂u}{∂ν}|^2\psi_1$ and hence in particular $\f{∂u}{∂ν}(x_0,t_0)=0$, which is a contradiction. 
We conclude that $\f{∂u}{∂ν}\geq 0$ on $\dOm\times(0,\infty)\cap\set{u>0}$ and hence, by continuity of $\f{∂u}{∂ν}$ and density of this set that $\f{∂u}{∂ν}\geq 0$ on $\dOm\times(0,\infty)$. 

Finally integrating \eqref{eq:utineq} over $\Om\times(0,t)$ and taking \eqref{12.1} into consideration, we see that 
\[
 \io u_0 \geq \io u(\cdot,t)\geq \io u_0 + \int_0^t\io Δu - χ\int_0^t\io \na\cdot\kl{\f uv\na v} =\io u_0+\int_0^t\intdom \f{∂u}{∂ν}
\]
by Gauss' theorem and $\f{∂v}{∂ν}=0$, which firstly shows that $\f{∂u}{∂ν}=0$ on $\dOm\times(0,∞)$ and secondly that \eqref{eq:utineq} actually is an equality. 
\qed
\mysection{Global smooth solutions to approximate problems}\label{sect3}
Now in order to approximate solutions by means of a convenient regularization of (\ref{01})-(\ref{02}),
for $\eps\in (0,1)$ we consider
\be{0eps}
	\left\{ \begin{array}{ll}
	u_{\eps t}=\Delta\ueps - \chi\nabla \cdot \Big(\frac{\ueps}{\veps} \nabla\veps\Big),
	\qquad & x\in \Omega, \ t>0, \\[1mm]
	v_{\eps t}=\Delta\veps-\veps+\frac{\ueps}{1+ε\ueps},
	\qquad & x\in \Omega, \ t>0, \\[1mm]
	\frac{\partial\ueps}{\partial\nu}=\frac{\partial\veps}{\partial\nu}=0,
	\qquad & x\in \pO, \ t>0, \\[1mm]
	\ueps(x,0)=u_0(x), \quad \veps(x,0)=v_0(x),
	\qquad & x\in \Omega,
	\end{array} \right.
\ee
and then first obtain the following.
\begin{lem}\label{lem15}
  For all $\eps\in (0,1)$, the problem (\ref{0eps}) admits a global classical solution 
  $(\ueps,\veps) \in (C^0(\bom\times [0,\infty))\cap C^{2,1}(\bom\times (0,\infty)))^2$ for which $\ueps>0$ in
  $\bom\times (0,\infty)$ and $\veps>0$ in $\bom\times [0,\infty)$.
\end{lem}
\proof
 The local existence of a solution can be obtained in a standard manner (cf. \cite[Lemma 3.1]{BBTW} for a related general setting). Boundedness of the source term $\f{\ue}{1+ε\ue}\leq \f1{ε}$ of the second equation in \eqref{0eps} translates into a bound on $t\mapsto\norm[W^{1,∞}(\Om)]{\ve(\cdot,t)}$. Combined with the strict positivity property of $\ve$ on $\Ombar\times(0,T)$ for any finite $T$ -- to be made more precise in Lemma \ref{lem2} below -- this serves to provide a uniform bound on $\ue$ on $\Ombar\times(0,T)$, in light of the extensibility criterion \cite[(3.3)]{BBTW} thus ensuring global existence of the solution. Positivity of $\ue$ follows from a classical strong maximum principle.
\qed
These approximate solutions clearly preserve mass:
\begin{lem}\label{lem1}
  Let $\eps\in (0,1)$. Then
  \be{mass}
	\io \ueps(\cdot,t) = \io u_0
	\qquad \mbox{for all } t>0.
  \ee
\end{lem}
\proof
  This directly results on integrating the first equation in (\ref{0eps}).
\qed
Moreover, the assumed positivity of $v_0$ enables us to control $\veps$ from below at least locally in time:
\begin{lem}\label{lem2}
  For each $\eps\in (0,1)$, we have
  \bas
	\veps(x,t)\ge \bigg( \inf_{x\in\Omega} v_0(x)\bigg) \cdot e^{-t}
	\qquad \mbox{for all $x\in\Omega$ and } t>0.
  \eas
\end{lem}
\proof
  As $\ueps$ is nonnegative, this is a straightforward consequence of a comparison argument applied to the second
  equation in (\ref{0eps}).
\qed
By means of well-known smoothing estimates of the heat semigroup, the mass conservation property (\ref{mass})
readily implies some basic regularity features of the second component.
\begin{lem}\label{lem3}
  Let $r \ge 1$ and $s\ge 1$ be such that $r<\frac{n}{n-2}$ and $s<\frac{n}{n-1}$. Then
  there exists $C>0$ such that for each $ε\in(0,1)$,
  \be{3.1}
	\io \veps^r(\cdot,t) \le C
	\qquad \mbox{for all } t>0
  \ee
  and
  \be{3.2}
	\io |\nabla\veps(\cdot,t)|^s \le C
	\qquad \mbox{for all } t>0.
  \ee
\end{lem}
\proof
 The representation of $\ve$ as 
\[
 \ve(\cdot,t)=e^{t(\Delta-1)} v_0 + \int_0^t e^{(t-s)(Δ-1)} \f{\ue(\cdot,s)}{1+ε\ue(\cdot,s)} ds
\]
makes it possible to apply well-known estimates for the Neumann heat-semigroup (cf. \cite[Lemma 1.3]{win_aggregationvs}), which provide positive constants $c_1, c_2, c_3$ and $c_4$ such that
\[
 \norm[\Lom r]{\ve(\cdot,t)} \leq c_1 \norm[\Lom r]{v_0} + c_2\int_0^t (1+(t-s)^{-\f n2(1-\f1r)} e^{-(t-s)}
	\Big\|\f{\ue(\cdot,s)}{1+ε\ue(\cdot,s)}\Big\|_{L^1(\Omega)}  ds
\]
and 
\[
 \norm[\Lom s]{\na \ve(\cdot,t)} \leq c_3\norm[W^{1,\infty}(\Om)]{v_0}+c_4\int_0^t (1+(t-s)^{-\f12-\f n2(1-\f 1s)}) e^{-(t-s)}\Big\|\f{\ue(\cdot,s)}{1+ε\ue(\cdot,s)}\Big\|_{L^1(\Omega)} ds
\]
for all $t>0$, so that the estimate $\norm[\Lom1]{\f{\ue}{1+ε\ue}}\leq \norm[\Lom1]{\ue(\cdot,t)}=\norm[\Lom1]{u_0}$ due to Lemma \ref{lem1} and finiteness of $\intninf (1+τ^{-\f n2(1-\f1r)}) e^{-τ} dτ$ and $\intninf (1+τ^{-\f12 -\f n2(1-\f1s)}) e^{-τ} dτ$ due to the conditions on $r$ and $s$ prove the lemma.
\qed
\mysection{A fundamental identity and first consequences thereof}\label{sect4}
Let us next formulate an identity which apparently reflects a fundamental structural propety of (\ref{01})-(\ref{02}),
as already used in a slightly modified form and for more restricted choices of $\chi$ in \cite{win_MMAS}.
In Lemma \ref{lem18} applied to $\varphi\equiv 1$, this will serve as a source for some essential a priori estimates
for (\ref{0eps}), whereas in \ref{lem77} we will make use of the freedom to choose widely arbitrary test functions
here in order to verify (\ref{9.3}) for the limit couple $(u,v)$ to be constructed in Lemma \ref{lem6}.
\begin{lem}\label{lem16}
  Let $p\in (0,1)$ and $q\in (0,1)$, and assume that $T>0$ and that $\varphi\in C^\infty(\bom\times [0,T])$
  is such that $\frac{\partial\varphi}{\partial\nu}=0$ on $\pO\times (0,T)$.
  Then 
  \bea{16.1}
	& & \hspace*{-20mm}
	-\int_0^T \io \ueps^p \veps^q \varphi_t
	+ \io \ueps^p(\cdot,T)\veps(\cdot,T) \varphi(\cdot,T)
	- \io u_0^p v_0^q \varphi(\cdot,0) \nn\\
	&=& \frac{4(1-p)q-4q^2-p(1-p)^2χ^2}{pq(pχ+1-q)} \int_0^T \io 
	\veps^q |\nabla \ueps^\frac{p}{2}|^2 \varphi \nn\\
	& &  + \frac{4(pχ+1-q)}{q} \int_0^T\io \bigg| \ueps^{\frac p2}\nabla \ve^{\frac q2} - \frac{(1-p)χ+2q}{2(pχ+1-q)}\ve^{\frac q2}\nabla \ue^{\frac p2}\bigg|^2φ\nn\\
	& & -\frac{2pχ}q \int_0^T \io 
	\ueps^\frac{p}{2} \veps^q \nabla\ueps^\frac{p}{2} \cdot\nabla\varphi \nn\\
	& & + \int_0^T \io \Big(1-\frac{p\chi}{q} \Big)\ueps^p \veps^q \Delta\varphi \nn\\
	& & - q \int_0^T \io \ueps^p \veps^q \varphi
	+ q \int_0^T \io \frac{\ueps^{p+1}\veps^{q-1}}{1+ε\ueps}  \varphi
	\qquad \mbox{for all } \eps\in (0,1).
  \eea
\end{lem}
\proof
  Using (\ref{0eps}), we compute
  \bea{16.2}
	\io \partial_t (\ueps^p \veps^q) \cdot \varphi
	&=& - p \io \nabla \Big(\ueps^{p-1} \veps^q \varphi\Big) \cdot 
		\Big(\nabla\ueps - \chi \frac{\ueps}{\veps}\nabla\veps\Big) \nn\\
	& & - q\io \nabla\Big(\ueps^p \veps^{q-1} \varphi\Big)\cdot\nabla\veps
	- q \io \ueps^p \veps^q \varphi
	+ q\io \frac{\ueps^{p+1}\veps^{q-1}}{1+ε\ueps} \varphi \nn\\[1mm]
	&=& p(1-p) \io \ueps^{p-2}|\nabla \ueps|^2\veps^qφ \nn\\
	& & - p(1-p) \chi \io \ueps^{p-1}\veps^{q-1} (\nabla\ueps\cdot\nabla\veps)\varphi 
	- 2pq \io \ueps^{p-1} \veps^{q-1} (\nabla\ueps\cdot\nabla\veps) \varphi \nn\\
	& & +pq\chi \io \ueps^p \veps^{q-2} |\nabla\veps|^2 \varphi
	+ q(1-q)\io \ueps^p \veps^{q-2} |\nabla\veps|^2 \varphi\nn\\
	& & - p \io \ueps^{p-1} \veps^q \nabla\ueps\cdot\nabla\varphi
	+ p\chi \io \ueps^p \veps^{q-1} \nabla\veps\cdot\nabla\varphi \nn\\
	& & -q \io \ueps^p \veps^{q-1} \nabla\veps\cdot\nabla\varphi \nn\\
	& & - q\io \ueps^p \veps^q \varphi
	+ q\io \frac{\ueps^{p+1} \veps^{q-1}}{1+ε\ueps} \varphi \nn\\[1mm]
	&=& \frac{4(1-p)}{p} \io \veps^q |\nabla\ueps^\frac{p}{2}|^2 \varphi \nn\\
	& & -\frac{4(χ-pχ+2q)}{q} \io \ueps^\frac{p}{2} \veps^\frac{q}{2} 
		(\nabla\ueps^\frac{p}{2}\cdot \nabla\veps^\frac{q}{2})\varphi\nn\\
	& & + \frac{4(p\chi+1-q)}{q} \io \ueps^p |\nabla\veps^\frac{q}{2}|^2 \varphi
	 -2\io \ueps^\frac{p}{2} \veps^q \nabla\ueps^\frac{p}{2} \cdot\nabla\varphi \nn\\
	& & + \bigg(\frac{p\chi}{q}-1\bigg) \io \ueps^p \nabla\veps^q \cdot\nabla\varphi\nn\\
	& & - q \io \ueps^p \veps^q \varphi
	+ q\io \frac{\ueps^{p+1}\veps^{q-1}}{1+ε\ueps} \varphi 
	\qquad \mbox{for all } t>0.
  \eea
  Here the same straighforward rearrangement as in \eqref{eq:rearrange} in the first three integrands on the right shows that
  \begin{align}\label{16.3}
&p(1-p) \io \ueps^{p-2}|\nabla \ueps|^2\veps^qφ 
 - p(1-p) \chi \io \ueps^{p-1}\veps^{q-1} (\nabla\ueps\cdot\nabla\veps)\varphi 
	- 2pq \io \ueps^{p-1} \veps^{q-1} (\nabla\ueps\cdot\nabla\veps) \varphi \nn\\ 
 &= \frac{4(pχ+1-q)}{q} \io \bigg|\ue^{\f p2}\na v^{\f q2}- \frac{(1-p)χ+2q}{2(pχ+1-q)}\ve^{\f q2}\na \ue^{\f p2}\bigg|^2 +\frac{4(1-p)q-4q^2-p(1-p)^2χ^2}{pq(pχ+1-q)}\io \ve^q|\na \ue^{\f p2}|^2.
  \end{align}
  Moreover, in the second of the two summands in (\ref{16.2}) which contain $\nabla\varphi$, we once more integrate by parts to 
  see that
 \begin{align*}
  -2\io \ue^{\f p2}\ve ^q\na\ue^{\f p2}\na φ+\Big(\frac{pχ}q-1\Big)\io \ue^p\na \ve ^q\na φ = \Big(1-\frac {pχ}q\Big)\io \ue ^p\ve^q\Delta φ - \frac{2pχ}q \io \ve ^q\ue ^{\f p2}\na\ue^{\f p2}\na φ
 \end{align*}
%
%
  thanks to the assumption that $\frac{\partial\varphi}{\partial\nu}=0$ on $\pO\times (0,T)$.
  Combining this with (\ref{16.2}) and (\ref{16.3}) establishes (\ref{16.1}).
\qed
An elementary but crucial observation now identifies a condition on the relationship between the exponents $p$ and 
$q$ which ensure positivity of the coefficient appearing in the first summand on the right-hand side in (\ref{16.1}).
\begin{lem}\label{lem166}
  Given $\chi>0$ and $p\in (0,1)$ such that $p<\frac{1}{\chi^2}$, let $\qp\in (0,1)$ and $\qm\in (0,\qp)$ be defined by
  \be{qpm}
	q_\pm(p):=\frac{1-p}{2} \cdot \Big(1\pm \sqrt{1-p\chi^2}\Big).
  \ee
  Then for any choice of $q\in (\qm,\qp)$,
  \be{166.1}
	\frac{4(1-p)q - 4q^2 - p(1-p)^2 \chi^2}{pq(p\chi+1-q)} >0
	\qquad \mbox{in } \Omega\times (0,\infty).
  \ee
\end{lem}
\proof
  We leave the straightforward proof to the reader.
\qed
%
%
%
%
%
As a consequence, for $p$ and $q$ as in Lemma \ref{lem166} we can readily derive the following
from Lemma \ref{lem16} when combined with the pointwise lower estimate for $\veps$ in Lemma \ref{lem2}.
\begin{lem}\label{lem18}
  Let $p\in (0,1)$ be such that $p<\frac{1}{\chi^2}$, and let $q\in (\qm,\qp)$ with $q_\pm(p)$ taken from (\ref{qpm}).
  Then for each $T>0$ there exists $C(p,q,T)>0$ fulfilling
  \be{18.1}
	\int_0^T \io \veps^q |\nabla\ueps^\frac{p}{2}|^2 \le C(p,q,T)
  \ee
  and
  \be{18.11}
	\int_0^T \io |\nabla\ueps^\frac{p}{2}|^2 \le C(p,q,T)
  \ee
  as well as
  \be{18.2}
 \frac{4(pχ+1-q)}{q} \int_0^T\io \bigg| \ueps^{\frac p2}\nabla \ve^{\frac q2} - \frac{(1-p)χ+2q}{2(pχ+1-q)}\ve^{\frac q2}\nabla \ue^{\frac p2}\bigg|^2
	\le C(p,q,T)
  \ee
  and
  \be{18.3}
	\int_0^T \io \ueps^{p+1} \veps^{q-1} \le C(p,q,T)
  \ee
  for all $\eps\in (0,1)$.
\end{lem}
\proof
  According to Lemma \ref{lem166}, our assumption $q\in (\qm,\qp)$ ensures that 
\[
 c_1:=\frac{4(1-p)q - 4q^2 - p(1-p)^2 \chi^2}{pq(p\chi+1-q)} 
\]
is positive. Moreover, Lemma \ref{lem2} along with (\ref{init}) says that given $T>0$ we can find
  $c_2>0$ such that
  \begin{align}\label{eq:vgeqc2}
	\veps(x,t) \ge c_2
	\qquad \mbox{for all $x\in\Omega$ and } t\in (0,T)
  \end{align}
  whenever $\eps\in (0,1)$, 
%
  whence applying Lemma \ref{lem16} to $\varphi\equiv 1$ shows that
  \bea{18.4}
	& & \hspace*{-10mm}
	c_1 c_2^q \int_0^T \io |\nabla\ueps^\frac{p}{2}|^2 \nn\\
	& & + \frac{4(pχ+1-q)}{q} \int_0^T\io \bigg| \ueps^{\frac p2}\nabla \ve^{\frac q2} - \frac{(1-p)χ+2q}{2(pχ+1-q)}\ve^{\frac q2}\nabla \ue^{\frac p2}\bigg|^2\nn\\
%
%
	& & +q \int_0^T \io \ueps^{p+1} \veps^{q-1} \nn\\
	&\le&  \io \ueps^p(\cdot,T)\veps^q(\cdot,T)
	- \io u_0^p v_0^q
	+ q\int_0^T \io \ueps^p \veps^q
	\qquad \mbox{for all } \eps\in (0,1).
  \eea
  Now by the H\"older inequality,
  \bas
	\io \ueps^p \veps^q \le \kkl{\io \ueps }^p \cdot \kkl{ \io \veps^\frac{q}{1-p}}^{1-p}
	\qquad \mbox{for all } t>0,
  \eas
  so that since
  \bas
	\frac{q}{1-p} < \frac{\qp}{1-p} 
	= \frac{1+\sqrt{1-p\chi^2}}{2} <1<\frac{n}{n-2},
  \eas
  we may combine (\ref{mass}) with Lemma \ref{lem3} to find $c_3>0$ fulfilling
  \bas
	\io \ueps^p \veps^q \le c_3
	\qquad \mbox{for all } t>0
  \eas
  whenever $\eps\in (0,1)$.
  The estimates in (\ref{18.1}), (\ref{18.2}) and (\ref{18.3}) therefore result from (\ref{18.4}),
  whereupon (\ref{18.11}) is a consequence of (\ref{18.1}) and \eqref{eq:vgeqc2}.
\qed
\mysection{A further consequence: A bound for $\ueps$ in $L^r$ for some $r>1$}\label{sect5}
Now in view of the desired integrability feature in (\ref{l1}), a crucial step in our analysis will consist
in deriving a spatio-temporal equi-integrability property of $\ueps$.
This will result from bounds therefor in some reflexive $L^r$ spaces, to be obtained by an interpolation between
(\ref{18.3}) and (\ref{3.1}). 
The following statement identifies the minimal possible choice of an integrability exponent arising in the course
of this argument (cf.~(\ref{5.100}) below), and will thereby form the core of our requirement (\ref{14.1}) on $\chi$.
%
%
\begin{lem}\label{lem55}
  Let $\chi>0$, and for $p\in(0,\min\{1,\frac1{χ^2}\}$ let $q_\pm(p)$ be as in (\ref{qpm}). Then
  \be{55.1}
	\inf_{\begin{array}{c} 
	\scriptstyle p\in (0,1), \ p<\frac{1}{\chi^2} \\[-1mm] \scriptstyle q\in (\qm,\qp) 
	\end{array}} \frac{1-q}{p}
	= \left\{ \begin{array}{ll}
	1 \qquad & \mbox{if } \chi\le 1, \\[1mm]
	\chi \qquad & \mbox{if } \chi\in (1,2), \\[1mm]
	1+\frac{\chi^2}{4}
	\qquad & \mbox{if } \chi\ge 2.
	\end{array} \right.
  \ee
\end{lem}
\proof
  By an evident monotonicity property,
  \bea{55.2}
	\inf_{\begin{array}{c} 
	\scriptstyle p\in (0,1), \ p<\frac{1}{\chi^2} \\[-1mm] \scriptstyle q\in (\qm,\qp) 
	\end{array}} \frac{1-q}{p}
	&=& \inf_{\scriptstyle p\in (0,1), \ p<\frac{1}{\chi^2}} 
	\frac{1-\qp}{p} \nn\\
	&=& \inf_{\scriptstyle p\in (0,1), \ p<\frac{1}{\chi^2}} 
	\frac{1-\frac{1-p}{2}(1+\sqrt{1-p\chi^2})}{p} \nn\\
	&=& \inf_{\scriptstyle p\in (0,1), \ p<\frac{1}{\chi^2}} 
	\frac{1+ p - (1-p)\sqrt{1-p\chi^2})}{2p}
	=: I(\chi)
  \eea
  for any $\chi>0$.
  Since $(1-p)\sqrt{1-p\chi^2}<1-p$ and thus $\frac{1+p-(1-p)\sqrt{1-p\chi^2}}{2p}>1$ for all
  $p\in (0,\min\{1,\frac{1}{\chi^2}\})$, and since on the other hand for $\chi\le 1$ we have
  \bas
	I(\chi) \le \liminf_{p\nearrow 1} \frac{1+p-(1-p)\sqrt{1-p\chi^2}}{2p}=1,
  \eas
  this firstly implies that $I(\chi)=1$ for any such $\chi$.\\
  In the case $\chi>1$, having in mind the substitution $\xi=\sqrt{1-p\chi^2}$ in (\ref{55.2}), we note that
  \bas
	\rho(\xi)
	&:=& \frac{1+\frac{1-\xi^2}{\chi^2} - (1-\frac{1-\xi^2}{\chi^2})\cdot\xi}{2\cdot \frac{1-\xi^2}{\chi^2}} \\
	&=& \frac{1}{2} \cdot \Big\{ \frac{\chi^2}{1-\xi^2} \cdot (1-\xi) + 1 + \xi \Big\} \\
	&=& \frac{1}{2} \cdot \Big\{ \frac{\chi^2}{1+\xi} + 1 + \xi \Big\},
	\qquad \xi\in [0,1),
  \eas
  satisfies
  \bas
	\rho'(\xi)=-\frac{\chi^2}{2(1+\xi)^2} + \frac{1}{2}
	\qquad \mbox{for all } \xi\in (0,1),
  \eas
  so that $\rho'$ attains a zero at $\xi=\chi-1\in (0,1)$ if and only if $\chi\in (1,2)$, while
  $\chi' \le 0$ throughout $(0,1)$ if $\chi\ge 2$.
  Therefore, $\inf_{\xi\in [0,1)} \rho(\xi)=\rho(\chi-1)=\chi$ if $\chi\in (1,2)$,
  whereas $\inf_{\xi\in [0,1)} \rho(\xi)=\lim_{\xi\nearrow 1} \rho(\xi)=1+\frac{\chi^2}{4}$ if $\chi\ge 2$.
  In conjunction with (\ref{55.2}), these observations verify (\ref{55.1}).
\qed
Now under the assumptions on $\chi$ from Theorem \ref{theo14}, the announced interpolation argument indeed
bears fruit of the desired flavour.
\begin{lem}\label{lem5}
  Suppose that $\chi>0$ is such that (\ref{14.1}) holds.
  Then there exists $r>1$ such that for any $T>0$ one can find $C(T)>0$ with the property that 
  \be{5.11}
	\int_0^T \io \ueps^r \le C(T)
	\qquad \mbox{for all } \eps\in (0,1).
  \ee
\end{lem}
\proof
  As a consequence of Lemma \ref{lem55}, our assumption on $\chi$ warrants that we can pick 
  $p\in (0,\min\{1,\frac{1}{\chi^2}\})$ and $q\in (\qm,\qp)$ such that
  \be{5.2}
	\frac{1-q}{p} <\frac{n}{n-2}.
  \ee
  Indeed, if $n=2$ this is obvious, while if $n\ge 4$ this is immediate from (\ref{55.1}), because then due to the 
  fact that $\frac{n}{n-2}\le 2$, the hypothesis (\ref{14.1}) in particular requires that $\chi<2$, so that in both 
  cases $\chi\le 1$ and $\chi>1$, (\ref{55.1}) shows that the assumption $\chi<\frac{n}{n-2}$ implies that
  \bas
	\inf_{\begin{array}{c} 
	\scriptstyle p\in (0,1), \ p<\frac{1}{\chi^2} \\[-1mm] \scriptstyle q\in (\qm,\qp) 
	\end{array}} \frac{1-q}{p}
	= \max\{1,\chi\}<\frac{n}{n-2}.
  \eas
  If $n=3$, in the case $\chi<2$ we similarly obtain that
  \bas
	\inf_{\begin{array}{c} 
	\scriptstyle p\in (0,1), \ p<\frac{1}{\chi^2} \\[-1mm] \scriptstyle q\in (\qm,\qp) 
	\end{array}} \frac{1-q}{p}
	= \max\{1,\chi\}<2<3=\frac{n}{n-2},
  \eas
  whereas when $\chi\ge 2$ we use our restriction $\chi<\sqrt{8}$ to infer from (\ref{55.1}) that
  \bas
	\inf_{\begin{array}{c} 
	\scriptstyle p\in (0,1), \ p<\frac{1}{\chi^2} \\[-1mm] \scriptstyle q\in (\qm,\qp) 
	\end{array}} \frac{1-q}{p}
	= 1+\frac{\chi^2}{4}<3
  \eas
  and that thus (\ref{5.2}) can be achieved also in this case.\abs
  Henceforth keeping $p$ and $q$ fixed such that (\ref{5.2}) holds, e.g.~by means of a continuity argument
  we can pick $r>1$ sufficiently close to $1$ such that still $p+1-r>0$ and
  \be{5.3}
	\frac{(1-q)r}{p+1-r} < \frac{n}{n-2}.
  \ee
  Then using Young's inequality, for $T>0$ we can estimate
  \bea{5.100}
	\int_0^T \io \ueps^r
	&=& \int_0^T \io \Big(\ueps^{p+1} \veps^{q-1} \Big)^\frac{r}{p+1} \cdot \veps^\frac{(1-q)r}{p+1} \nn\\
	&\le& \int_0^T \io \ueps^{p+1} \veps^{q-1} 
	+ \int_0^T \io \veps^\frac{(1-q)r}{p+1-r}
	\qquad \mbox{for all } \eps\in (0,1),
  \eea
  so that (\ref{5.11}) results on using (\ref{18.3}) and applying (\ref{3.1}) together with (\ref{5.3}).
\qed
\mysection{A weighted $L^2$ bound for $\nabla\veps$}\label{sect6}
In order to complement (\ref{18.11}) by an analogous $L^2$ estimate for $\nabla \veps$ merely involving $\veps$
but not $\ueps$ as a weight function, 
independently from the above we apply a standard testing technique to the second equation in (\ref{0eps}) with
the following outcome.
\begin{lem}\label{lem17}
  For all $q\in (0,1)$ and any $T>0$ one can find $C(T)>0$ such that
  \be{17.1}
	\int_0^T \io |\nabla\veps^\frac{q}{2}|^2 \le C(T)
	\qquad \mbox{for all } \eps\in (0,1).
  \ee
\end{lem}
\proof
  Thanks to the positivity of $\veps$, we may use $\veps^{q-1}$ as a test function in the second equation in (\ref{0eps})
  to see that
  \bea{17.2}
	\frac{1}{q} \frac{d}{dt} \io \veps^q
	&=& (1-q)\io \veps^{q-2} |\nabla \veps|^2 - \io \veps^q + \io \frac{\ueps}{1+ε\ue} \veps^{q-1} \nn\\
	&\ge& (1-q)\io \veps^{q-2} |\nabla \veps|^2 - \io \veps^q
	\qquad \mbox{for all } t>0,
  \eea
  where according to Lemma \ref{lem3} and the fact that $q<1<\frac{n}{n-2}$, we can find $c_1>0$ such that
  \bas
	\io \veps^q \le c_1
	\qquad \mbox{for all } t>0.
  \eas
  On integration, we thus obtain from (\ref{17.2}) that
  \bas
	\frac{4(1-q)}{q^2} \int_0^T \io |\nabla\veps^\frac{q}{2}|^2
	&=& (1-q) \int_0^T \io \veps^{q-2}|\nabla\veps|^2 \\
	&\le& \frac{1}{q} \io \veps^q(\cdot,T) + \int_0^T \io \veps^q \\
	&\le& \frac{c_1}{q} + Tc_1
  \eas
  for all $\eps\in (0,1)$.
\qed
\mysection{Time regularity}\label{sect7}
As a final preparation for our limit procedure, we establish some regularity features of the time derivatives
in (\ref{0eps}), beginning with a conveniently transformed version of the first solution component.
\begin{lem}\label{lem34}
  Assume (\ref{14.1}), and let $p\in (0,1)$ be such that $p<\frac{1}{\chi^2}$. Then for all $T>0$ there exists $C(T)>0$
  such that
  \be{34.1}
	\int_0^T \Big\|\partial_t \Big(\ueps(\cdot,t)+1\Big)^\frac{p}{2}\Big\|_{(W_0^{1,\infty}(\Omega))^\star} dt
	\le C(T)
	\qquad \mbox{for all } \eps\in (0,1).
  \ee
\end{lem}
\proof
  We fix $\psi\in C_0^\infty(\Omega)$ such that $\|\psi\|_{W^{1,\infty}(\Omega)} \le 1$ and use (\ref{0eps}) and
  Young's inequality as well as the trivial estimates $\ue\le (\ue+1)$ and $(\ue+1)^{-\f p2}\le 1$ to see that
  \bas
	\bigg| \io \partial_t (\ueps+1)^\frac{p}{2}\psi \bigg|
	&=& \bigg| \frac{p(2-p)}{4} \io (\ueps+1)^\frac{p-4}{2} |\nabla\ueps|^2 \psi
	- \frac{p}{2} \io (\ueps+1)^\frac{p-2}{2} \nabla\ueps\cdot\nabla\psi \\
	& & - \frac{p(2-p)\chi}{4} \io \frac{\ueps(\ueps+1)^\frac{p-4}{2}}{\veps} 
		(\nabla\ueps\cdot\nabla\veps) \psi 
	+ \frac{p\chi}{2} \io \frac{\ueps (\ueps+1)^\frac{p-2}{2}}{\veps} \nabla\veps\cdot\nabla\psi 
	\bigg| \\[1mm]
	&\le& \frac{p(2-p)}{4} \io (\ueps+1)^{\frac{p-4}{2}} |\nabla\ueps|^2
	+ \frac{p}{2} \io (\ueps+1)^\frac{p-2}{2} |\nabla\ueps| \\
	& & + \frac{p(2-p)\chi}{4} \io (\ueps+1)^\frac{p-2}{2} |\nabla\ueps| \cdot \frac{|\nabla\veps|}{\veps}
	+ \frac{p\chi}{2} \io (\ueps+1)^\frac{p}{2} \frac{|\nabla\veps|}{\veps} \\[1mm]
	&\le& \frac{p(2-p)}{4} \io (\ueps+1)^{p-2} |\nabla\ueps|^2 
	+ \frac{p}{4} \io (\ueps+1)^{p-2} |\nabla\ueps|^2 
	+ \frac{p|\Omega|}{4} \\
	& & + \frac{p(2-p)\chi}{8} \io (\ueps+1)^{p-2} |\nabla\ueps|^2
	+ \frac{p(2-p)\chi}{8} \io \frac{|\nabla\veps|^2}{\veps^2} \\
	& & + \frac{p\chi}{4} \io (\ueps+1)^p
	+ \frac{p\chi}{4} \io \frac{|\nabla\veps|^2}{\veps^2}
	\qquad \mbox{for all $t>0$ and } \eps\in (0,1).
  \eas
  Since Lemma \ref{lem2} provides $c_1>0$ such that $\veps\ge c_1$ in $\Omega\times (0,T)$ and hence
  \bas
	\io \frac{|\nabla\veps|^2}{\veps^2} \le \frac{4}{c_1} \io |\nabla\veps^\frac{1}{2}|^2
	\qquad \mbox{for all $t\in (0,T)$ and } \eps\in (0,1),
  \eas
  and since
  \bas
	\io (\ueps+1)^p \le \io u_0+|\Omega|
	\qquad \mbox{for all $t>0$ and } \eps\in (0,1)
  \eas
  by (\ref{mass}), we thus infer that there exists $c_2>0$ such that for all $\eps\in (0,1)$,
  \bas
	\Big\|\partial_t \Big(\ueps(\cdot,t)+1\Big)^\frac{p}{2}\Big\|_{(W_0^{1,\infty}(\Omega))^\star}
	&=& \sup_{ \begin{array}{c}
	\scriptstyle \psi\in C_0^\infty(\Omega) \\[-1mm] \scriptstyle \|\psi\|_{W^{1,\infty}(\Omega)} \le 1
	\end{array}} 
	\bigg| \io \partial_t \Big(\ueps(\cdot,t)+1\Big)^\frac{p}{2} \psi \bigg| \\
	&\le& c_2 \cdot \kkl{ \io |\nabla\ueps^\frac{p}{2}|^2
	+ \io |\nabla\veps^\frac{1}{2}|^2 + 1 }
	\qquad \mbox{for all } t\in (0,T).
  \eas
  Thanks to the outcomes of Lemma \ref{lem18} and Lemma \ref{lem17}, an integration over $t\in (0,T)$ therefore
  yields (\ref{34.1}).
\qed
As for the second component, we can directly address the quantity $v_{\eps t}$.
\begin{lem}\label{lem35}
  Let $\chi>0$. Then there exists $C>0$ such that whenever $\eps\in (0,1)$,
  \be{35.1}
	\|v_{\eps t}(\cdot,t)\|_{(W_0^{1,\infty}(\Omega))^\star} dt
	\le C
	\qquad \mbox{for all } t>0.
  \ee
\end{lem}
\proof
  We again fix $\psi\in C_0^\infty(\Omega)$ fulfilling $\|\psi\|_{W^{1,\infty}(\Omega)}\le 1$, and using
  (\ref{0eps}) we find that
  \bas
	\bigg| \io v_{\eps t} \psi\bigg|
	&=& \bigg| - \io \nabla\veps\cdot\nabla \psi - \io \veps \psi + \io \f{\ueps}{1+ε\ue} \psi \bigg| \\
	&\le& \io |\nabla\veps| + \io \veps + \io \ueps
	\qquad \mbox{for all $t>0$ and } \eps\in (0,1).
  \eas
  Therefore,
  \begin{equation*}
	\|v_{\eps t}(\cdot,t)\|_{(W_0^{1,\infty}(\Omega))^\star}
	\le \sup_{\tau>0} \bigg\{ \io |\nabla\veps(\cdot,\tau)| + \io \veps (\cdot,\tau) + \io \ueps(\cdot,\tau)\bigg\}
  	\qquad \mbox{for all $t>0$ and } \eps\in (0,1),
  \end{equation*}
  so that (\ref{35.1}) results from Lemma \ref{lem3} and (\ref{mass}).
\qed
\mysection{Construction of limit functions. Proof of Theorem \ref{theo14}}\label{sect8}
Collecting the above estimates, by means of a straightforward extraction procedure we can pass to the limit
$\eps\searrow 0$ in the following sense.
\begin{lem}\label{lem6}
  Suppose that (\ref{14.1}) holds, and let $p\in (0,1)$ and $q\in (0,1)$ be such that $p<\frac{1}{\chi^2}$
  and $q\in (\qm,\qp)$. Then 
  there exist $(\eps_j)_{j\in\N}\subset (0,1)$ and functions $u$ and $v$ defined on $\Omega\times (0,\infty)$
  such that $\eps_j\searrow 0$ as $j\to\infty$, that $u\ge 0$ and $v>0$ a.e.~in $\Omega\times (0,\infty)$, and that
  \begin{align}
	&\qquad && \ueps\to u
	&&\mbox{in $L^1_{loc}(\bom\times [0,\infty))$ and a.e.~in } \Omega\times (0,\infty),
	\label{6.1} \\
	&&& \nabla\ueps^\frac{p}{2} \wto \nabla u^\frac{p}{2}
	&&\mbox{in } L^2_{loc}(\bom\times [0,\infty)),
	\label{6.4} \\
	&&& \veps \to v
	&&\mbox{in $L^1_{loc}(\bom\times [0,\infty))$ and a.e.~in } \Omega\times (0,\infty),
	\label{6.2} \\
	&&& \nabla\veps \wto \nabla v
	&&\mbox{in } L^1_{loc}(\bom\times [0,\infty))
	\qquad \mbox{and}
	\label{6.3} \\
	&&& \nabla\veps^\frac{q}{2} \wto \nabla v^\frac{q}{2}
	&&\mbox{in } L^2_{loc}(\bom\times [0,\infty))
	\label{6.5}
  \intertext{as $\eps=\eps_j\searrow 0$. Moreover,}
	&\qquad&&\ve^{\f q2}\na\ue^{\f p2} \wto v^{\f q2} \na u^{\f p2}
	&&\text{ in } L^2_{loc}(\Ombar\times[0,∞)),
	\label{eq:vq2naup2weakL2} \\
	&&& \ue^{\f p2}\na \ve^{\f q2}\wto u^{\f p2}\na v^{\f q2} 
	&&\text{ in } L^2_{loc}(\Ombar\times[0,∞)),
	\label{eq:up2navq2weakL2} \\
 	&&&  \ue^p\ve^q \to u^pv^q  
	&&\text{ in } L^1_{loc}(\Ombar\times[0,∞)),
	\label{eq:upvqL1} \\
	&&& \ue^{\f p2}\ve^{\f q2} \to u^{\f p2}v^{\f q2}  
	&&\text{ in } L^2_{loc}(\Ombar\times[0,∞)) \qquad \mbox{as well as}
	\label{eq:up2vq2L2} \\
	&&& \ue^{p+1}\ve^{q-1} \to u^{p+1}v^{q-1} 
	&&\text{ in } L^1_{loc}(\Ombar\times[0,∞))
	\label{eq:upp1vqm1L1}
  \end{align}
  as $\eps=\eps_j\searrow 0$, and 
  \be{6.6}
	\io u(\cdot,t) = \io u_0
	\qquad \mbox{for a.e.~} t>0.
  \ee
\end{lem}
\proof
  We fix $p\in (0,1)$ such that $p<\frac{1}{\chi^2}$ and combine Lemma \ref{lem18} with (\ref{mass}) and Lemma \ref{lem34}
  to see that
  \bas
	\Big((\ueps+1)^\frac{p}{2}\Big)_{\eps\in (0,1)}
	\quad \mbox{is bounded in } L^2_{loc}([0,\infty);W^{1,2}(\Omega))
  \eas
  and that
  \bas
	\Big( \partial_t (\ueps+1)^\frac{p}{2}\Big)_{\eps\in (0,1)}
	\quad \mbox{is bounded in } L^1_{loc}([0,\infty);(W_0^{1,\infty}(\Omega))^\star).
  \eas
  Apart from that, Lemma \ref{lem3} and Lemma \ref{lem35} show that there exists $r>1$ such that
  \bas
	(\veps)_{\eps\in (0,1)}
	\quad \mbox{is bounded in } L^r_{loc}([0,\infty);W^{1,r}(\Omega))
  \eas
  and
  \bas
	(v_{\eps t})_{\eps\in (0,1)}
	\quad \mbox{is bounded in } L^\infty((0,\infty);(W_0^{1,\infty}(\Omega))^\star).
  \eas
  Therefore, by means of two applications of an Aubin-Lions lemma \cite[Cor. 8.4]{simon} we can find $(\eps_j)_{j\in\N}\subset (0,1)$
  such that $\eps_j\searrow 0$ as $j\to\infty$, that $\ueps^\frac{p}{2} \to u^\frac{p}{2}$ in
  $L^2_{loc}(\bom\times [0,\infty))$ and a.e.~in $\Omega\times (0,\infty)$ 
  as $\eps=\eps_j\searrow 0$, and that
  (\ref{6.4}), (\ref{6.2}) and (\ref{6.3}) hold with some nonnegative functions $u$ and $v$ defined on
  $\Omega\times (0,\infty)$.
  Since thus also $\ueps\to u$ a.e.~in $\Omega\times (0,\infty)$ as $\eps=\eps_j\searrow 0$,
  making use of the equi-integrability property implied by Lemma \ref{lem5} we may invoke the Vitali convergence theorem
  to infer that in fact also (\ref{6.1}) holds, whereupon (\ref{6.6}) becomes a consequence of (\ref{mass}).
  The additional convergence statement in (\ref{6.5}) finally results from Lemma \ref{lem17} and (\ref{6.2})
  in a straightforward manner.\\
Because $q<1$, \eqref{6.2} and Vitali's convergence theorem imply that $\ve^{\f q2}\to v^{\f q2}$ in $\LomT2$. Combining this with \eqref{6.4}, we obtain $\ve^{\f q2}\na\ue^{\f p2}\wto v^{\f q2} \na u^{\f p2}$ in $\LomT 1$ and the finiteness of $\sup_{ε\in(0,1)} \norm[L^2(\Om\times(0,T))]{\ve^{\f q2}|\na \ue^{\f p2}|}$ according to \eqref{18.1} serves to prove \eqref{eq:vq2naup2weakL2}. Analogously invoking \eqref{6.1} and \eqref{6.5}, we infer 
\begin{equation}\label{eq:up2navq2weakL1} 
\ue^{\f p2}\na \ve^{\f q2}\wto u^{\f p2}\na v^{\f q2}\quad \text{ in } \LomT1.
\end{equation}
 According to \eqref{18.2}, for some $c_1\in ℝ$, $\ueps^{\frac p2}\nabla \ve^{\frac q2} - c_1\ve^{\frac q2}\nabla \ue^{\frac p2}$ is bounded in $\LomT2$ and hence, due to \eqref{eq:vq2naup2weakL2}, so is $\ue^{\f p2}\na \ve^{\f q2}$ so that we can sharpen \eqref{eq:up2navq2weakL1} into \eqref{eq:up2navq2weakL2}. 
 By Young's inequality, (\ref{mass}) and Lemma \ref{lem3} there exists $c_2>0$ such that
  \bas
	\int_0^T \io \Big| \ueps^p \veps^q \Big|^\frac{1}{p+q}
	\le \int_0^T \io \ueps + \int_0^T \io \veps \le c_2
	\qquad \mbox{for all } \eps\in (0,1)
  \eas
  and hence Vitali's convergence theorem proves \eqref{eq:upvqL1}, because $p+q<p+\qp<1$. This also immediately implies \eqref{eq:up2vq2L2}. 
 For proving \eqref{eq:upp1vqm1L1} we use the continuity of $q_{\pm}(\cdot)$ to pick $\ptilde>p$ such that $\ptilde<1$, $\ptilde<\f1{χ^2}$ and 
\[
 \qtilde:=1+\f{\ptilde+1}{p+1}(q-1)\in (q_-(\ptilde),q_+(\ptilde)). 
\]
Then by Lemma \ref{lem3} we can find $c_3>0$ such that 
\[
 \intnT\io (\ue^{p+1}\ve^{q-1})^{\f {\ptilde+1}{p+1}} = \intnT\io \ue^{\ptilde+1}\ve^{\qtilde-1} \leq c_3 
\]
and the Vitali convergence theorem shows \eqref{eq:upp1vqm1L1}. 
\qed
Our next aim is to make sure that the functions $u$ and $v$ we have just constructed form a generalized solution of \eqref{01}-\eqref{02}. We begin with the second equation. 
\begin{lem}\label{lem7}
  If (\ref{14.1}) holds, then the pair $(u,v)$ obtained in Lemma \ref{lem6} is a global weak solution of (\ref{02})
  in the sense of Definition \ref{defi11}.
\end{lem}
\proof
  From (\ref{6.1}), (\ref{6.2}) and (\ref{6.3}) we immediately see that the regularity properties in (\ref{11.1})
  hold, and that moreover for arbitrary $\varphi\in C_0^\infty(\bom\times [0,\infty))$, in the identity
  \bas
	- \int_0^\infty \io  \veps \varphi_t - \io v_0 \varphi(\cdot,0)
	= - \int_0^\infty \io \nabla\veps\cdot\nabla\varphi
	-\int_0^\infty \io \veps\varphi
	+ \int_0^\infty \io \ueps\varphi,
  \eas
  valid for all $\eps\in (0,1)$ due to (\ref{0eps}), we may let $\eps=\eps_j\searrow 0$ in each integral separately
  to readily verify (\ref{11.2}).
\qed
Another important part of Definition \ref{defi12} are positivity requirements, which will be established 
in Lemma \ref{lem:uvpositive}. 
The following technical lemmas prepare the main argument therein, 
where we will derive a differential inequality for $\io \ln \ue$,
and where further exploting the latter will in particular require some
$\eps$-independent lower bound for this functional at some suitable initial value, 
despite the fact that \eqref{init} does not guarantee finiteness of $\io \ln u_0$. 
An appropriate replacement, to be provided by Lemma \ref{lem:iolnue}, is entailed by the comparison-type Lemma \ref{lem:diffineq} in combination with a differential inequality, the derivation of which rests on Lemma \ref{lem:poincareforlog}. 
\begin{lem}\label{lem:diffineq}
  Let $a>0$, $b>0$ and $T>0$, and let $y\colon(0,T)\to ℝ$ be a continuously differentiable function satisfying 
  \[
 	y'(t)\leq -ay^2(t)+b\qquad \text{for all } t\in(0,T) \text{ at which } y(t)>0. 
  \]
  Then 
  \[
 	y(t)\leq \sqrt{\f ba}\coth(\sqrt{ab}t) \qquad \text{for all } t\in(0,T).
  \]
\end{lem}
\proof
  Let $η>0$. Then $M_{η}:=\Big\{t\in(0,T): y(t)>\sqrt{\f ba}+η\Big\}$ (is either empty or) can be written as union of 
  its connected components, i.e. $M_{η}=\bigcup_{k\inℕ} I_k$ with disjoint open intervals $I_k$. If we consider any
  nonempty $I_k$ with $\inf I_k\neq 0$, by continuity $y(\inf I_k)=\sqrt{\f ba}+η$ and hence $y'(\inf I_k)=-2η-aη^2<0$, 
  contradicting the definition of $\inf I_k$ as infimum of a set where $y>\sqrt{\f ba}+η$. Hence there is $t_{η}\in[0,T)$ 
  such that $y\leq \sqrt{\f ba}+η$ in $(t_{η},T)$ and that $y>\sqrt{\f ba}+η$ in $(0,t_{η})$ so that $b-ay^2$ 
  is negative in $(0,t_{η})$ and for $t_0\in(0,t_{η})$ and $t\in (t_0,t_\eta)$ we find that 
  \[
 	t-t_0 = \int_{t_0}^t \f{y'(s)}{b-ay^2(s)} ds 
	=\f1{\sqrt{ab}}\int_{\sqrt{\f ab}y(t_0)}^{\sqrt{\f ab}y(t)} \f1{1-z^2} dz
	=\f1{\sqrt{ab}}\kkl{ \arcoth\Big(\sqrt{\f ab} y(t)\Big) - \arcoth\Big(\sqrt{\f ab}y(t_0)\Big)},
  \]
  leading to 
  \[
 	\sqrt{ab}(t-t_0) + \arcoth\Big(\sqrt{\f ab}y(t_0)\Big) \leq \arcoth\Big(\sqrt{\f ab} y(t)\Big)
  \]
  and hence to 
  \[
 	y(t) \leq \sqrt{\f ba} \coth\bigg(\sqrt{ab}(t-t_0) + \arcoth\Big(\sqrt{\f ab}y(t_0)\Big)\bigg) 
	\leq \sqrt{\f ba} \coth\Big(\sqrt{ab}(t-t_0)\Big).
  \]
  Using that $t_0\in(0,t)$ and $η>0$ were arbitrary, we conclude that 
  $y(t)\leq \max\Big\{\sqrt{\f ba}\coth(\sqrt{ab}t), \sqrt{\f ba}\Big\} = \sqrt{\f ba}\coth(\sqrt{ab}t)$.
\qed
The following statement essentially goes back to an observation made in \cite{taowin_persistenceofmass}.
\begin{lem}\label{lem:poincareforlog}
  Let $η>0$. Then there exists $C>0$ such that every positive function $φ\in C^1(\Ombar)$ fulfilling 
  \[
  	\Big|\set{x\in \Om; \; φ(x)>δ} \Big|>η
  \]
  for some $δ>0$ satisfies 
  \[
 	\io \f{|\na φ|^2}{φ^2} \geq C \cdot \kkl{\io \ln\f{δ}{φ}}^2 \quad \text{or}\quad \io \ln \f{δ}{φ}<0.
  \]
\end{lem}
\proof
  This directly follows from the inequality provided by \cite[Lemma 4.3]{taowin_persistenceofmass} upon squaring. 
  A requirement on convexity of the domain, as additionally made there in order to allow for an application of
  the Poincar\'e inequality from \cite[Cor 9.1.4]{jost_pde} in the proof, can actually be removed by replacing the
  latter with Lemma \ref{lem:poincare}.
\qed
We can now pass to our derivation of lower bounds for $\io \ln \ueps$ in the following form.
\begin{lem}\label{lem:iolnue}
  There exists $T>0$ such that for every $t\in(0,T)$, 
  \[
 	\inf_{ε\in(0,1)} \io \ln \ue (\cdot,t) > -∞.
  \]
\end{lem}
\proof
We let $M_{ε}(t):=\sup_{τ\in[0,t]} \norm[\Lom\infty]{\ue(\cdot,τ)}$ for $t\in (0,∞)$ and $\eps\in (0,1)$,
and pick $p>n$. From 
Lemma 1.3 iii) and Lemma 1.3 ii) in \cite{win_aggregationvs} we obtain $c_1>0$ and $c_2>0$ such that 
\bas
 \norm[\Lom\infty]{\na \ve(\cdot,t)}
&\leq& c_1\norm[W^{1,\infty}(\Om)]{v_0} + c_2\intnt (1+(t-s)^{-\f12})e^{-λ(t-s)}\norm[\Lom\infty]{\ue} ds  \\
&\leq& c_3 \cdot (1+M_\eps(t))
	\qquad \mbox{for all } t\in (0,∞),
\eas
where $c_3=\max\Big\{c_1 \|v_0\|_{W^{1,\infty}(\Omega)},c_2\intninf (1+τ^{-\f12})e^{-λτ}dτ\Big\}$. 
Invoking Lemma 1.3 (iv) from \cite{win_aggregationvs} 
we find $c_4>0$ such that 
\begin{align*}
 \norm[\Lom\infty]{u(\cdot,t)}
&\leq \norm[\Lom\infty]{u_0} + χc_4\intnt (1+(t-s)^{-\f12-\f n{2p}})\norm[\Lom p]{\f{\ue}{\ve}\na \ve} ds\\
&\leq \norm[\Lom\infty]{u_0} + \f{χc_4}{\inf v_0} e^t \cdot c_3(1+M_\eps(t)) \cdot \norm[\Lom 1]{u_0} (M_\eps(t))^{\f{p-1}p}
\intnt (1+τ^{-\f12-\f n{2p}}) dτ
\end{align*}
for all $t\in (0,∞)$, so that in conclusion we can find $c_5>0$ fulfilling
\[
 M_\eps(t)\leq \norm[\Lom\infty]{u_0} + c_5 (t+t^{\f12-\f n{2p}})e^t (1+M_\eps(t)) (M_\eps(t))^{\f {p-1}p} 
	\qquad \mbox{for all } t\in (0,∞),
\]
and hence by the fact that for all $a,b \in [0,∞)$, $γ\in(0,1)$
\[
 \sup\set{x\in[0,∞);\; x\leq a+bx^{γ}}\leq \f{a}{1-γ}+b^{\f1{1-γ}}
\]
we can achieve that 
\[
 M_\eps(t) \leq p\norm[\Lom\infty]{u_0} + \bigg(c_5(t+t^{\f12-\f n{2p}})e^t (1+M_\eps(t)) \bigg)^p
\]
If we let $T_{ε}:=\sup\Big\{t\in (0,∞): M_\eps(t)\leq p\norm[\Lom\infty]{u_0}+1\Big\}$, then certainly 
\[
T_{ε}>\min\bigg\{1,\f1{2 c_5 \cdot 2e (2+p\norm[\Lom\infty]{u_0})}\bigg\}=:T.
\]
 In conclusion, this means that for all $ε\in(0,1)$,
\[
 \norm[\Lom\infty]{\ue(\cdot,t)}\leq p\norm[\Lom\infty]{u_0}+1=:M \quad\text{and}\quad \norm[W^{1,\infty}(\Om)]{\na\ve(\cdot,t)}\leq c_3(1+M) \qquad \text{for all } t\in(0,T).
\]
In particular with $δ:=\f{1}{2|\Om|}\io u_0$ and $η:=\f1{2M}\io u_0$ we have $|\set{\ue(\cdot,t)\geq δ}|\geq η$ for every $t\in(0,T)$ and each $ε\in(0,1)$. From Lemma \ref{lem:poincareforlog} we hence obtain $c_6>0$ such that 
\begin{align*}
\frac{d}{dt}
\bigg(\io \ln \f{δ}{\ue(\cdot,t)}\bigg) &= -\io \f{|\na \ue(\cdot,t)|^2}{\ue^2(\cdot,t)} + χ\io \f1{\ue(\cdot,t)\ve(\cdot,t)} \na\ue(\cdot,t)\cdot\na\ve(\cdot,t)\\
&\leq -\f12\io \f{|\na \ue(\cdot,t)|^2}{\ue^2(\cdot,t)} + \f{χ^2}2\io \f{|\na \ve(\cdot,t)|^2}{\ve^2(\cdot,t)}\\ 
&\leq  -\f{c_6}2 \bigg(\io \ln\f{δ}{\ue(\cdot,t)}\bigg)^2 + \f{χ^2e^{2T}c_3^2(1+M)^2}{2(\inf v_0)^2}
\end{align*}
for every $t\in(0,T)$ at which $\io \ln\f{δ}{\ue(\cdot,t)}>0$. Lemma \ref{lem:diffineq} hence proves the claim. 
\qed
This enables us to verify the positivity requirements from Definition \ref{defi12} without any assumptions on the initial
data beyond (\ref{init}).
\begin{lem}\label{lem:uvpositive}
  The functions $u$, $v$ obtained in Lemma \ref{lem6} satisfy 
  $v>0$, $u>0$ a.e. in $\Om\times(0,\infty)$ and $u^pv^q>0$ a.e. on $\dOm\times(0,\infty)$.
\end{lem}
\proof
  According to Lemma \ref{lem2} and \eqref{6.2}, $\essinf_{x\in\Omega} v(x,t)>\inf v_0 e^{-t}$
  for any $t>0$, 
  and \eqref{6.4} and \eqref{6.5} together with \eqref{6.1} or \eqref{6.2}, respectively, show that $u$ and $v$ can be 
  evaluated on $\dOm\times(0,\infty)$ in the sense of traces. For the proof of the remaining positivity properties $u>0$ 
  a.e. in $\Om\times(0,\infty)$ and $u>0$ a.e. on $\dOm\times(0,\infty)$, we intend to prove 
  \begin{equation}\label{eq:lnuinw12}
 	\ln u \in L^2_{loc}((0,\infty);W^{1,2}(\Om)), 
  \end{equation}
  which entails $\ln u\in L^2_{loc}(\Ombar\times(0,\infty))$ and, due to the embedding 
  $W^{1,2}(\Om)\hookrightarrow L^2(\dOm)$, also $\ln u\in L^2_{loc}(\dOm)\times(0,\infty))$, proving positivity of $u$ 
  a.e. in the respective sets. By Lemmata \ref{lem2} and \ref{lem1},
  \begin{align*}
 	\frac{d}{dt} \bigg[-\io \ln\ue - χ^2\io \ln \ve\bigg] 
	&= -\io \f{|\na \ue|^2}{\ue^2} + χ\io \f{1}{\ue\ve}\na\ue\cdot\na \ve - χ^2 \io \f{|\na\ve|^2}{\ve^2}\\
	&\qquad +χ^2|\Om| - χ^2\io \f{\ue}{\ve(1+ε\ue)}\\
 	&\leq  -\f12\io \f{|\na \ue|^2}{\ue^2} - \f{χ^2}{2} \io \f{|\na \ve|^2}{\ve^2} + χ^2|\Om| 
	+ \f{χ^2}{\inf v_0} e^t \io u_0
  \end{align*}
  for any $t>0$. 
  Now picking $\tau>0$, according to Lemma \ref{lem:iolnue} we can find $\tau_0\in (0,\tau)$
  and $c_1>0$ such that
  \be{foo}
	\inf_{\eps\in(0,1)} \io \ln \ueps(\cdot,\tau_0) >-\infty. 
  \ee
  For any fixed $T>τ$ we then obtain 
  \begin{equation}\label{eq:minuslnu}
 	-\io \ln \ue(\cdot,t) + \f12 \int_{τ_0}^t\io |\na \ln \ue|^2 
	\leq -\io \ln \ue(\cdot,τ_0) +χ^2\io \ln\f{\ve(\cdot,t)}{\ve(\cdot,τ_0)} + χ^2|\Om| + \f{χ^2}{\inf v_0} e^T \io u_0 
  \end{equation}
  for $t\in(τ,T)$, where according to \eqref{foo} and by Lemma \ref{lem2} and \eqref{3.1} the right-hand side is bounded independently of $ε$. We invoke the Poincar\'e inequality to find $c_1>0$ such that 
$\norm[W^{1,2}(\Om)]{φ}\leq c_1(\norm[\Lom2]{\na φ}+\norm[\Lom1]{φ})$ for all $φ\in W^{1,2}(\Om)$, 
and since the elementary estimate $|\ln s|\leq 2s-\ln s$ valid for all $s>0$, Lemma \ref{lem1} and \eqref{eq:minuslnu} provide $c_2>0$ such that $\io |\ln \ue(\cdot,t)|\leq c_2\) for all $t\in(τ,T)$, 
  we conclude that for every $τ>0$ and $T>τ$ there exists $c_3>0$ such that for all $ε\in(0,1)$ we have
  \[
 	\norm[L^2((τ,T);W^{1,2}(\Om))]{\ln \ue} \leq c_3,
  \]
  which by a weak compactness argument immediately results in \eqref{eq:lnuinw12}.
\qed
We can now make sure that indeed the obtained pair $(u,v)$ has all the properties required in 
Definition \ref{defi9}.
\begin{lem}\label{lem77}
  Suppose that (\ref{14.1}) holds, and let $p\in (0,1)$ and $q\in (0,1)$ be such that $p<\frac{1}{\chi^2}$
  and $q\in (\qm,\qp)$.
  Then the functions $u$ and $v$ constructed in Lemma \ref{lem6} form a global weak $(p,q)$-supersolution
  of (\ref{01}) in the framework of Definition \ref{defi9}.
\end{lem}
\proof
  The regularity requirements in \eqref{9.1} are satisfied according to \eqref{eq:upvqL1} and \eqref{eq:upp1vqm1L1},  
  whereas those in \eqref{9.2} result from \eqref{eq:vq2naup2weakL2} and \eqref{eq:up2navq2weakL2}. Positivity a.e. has 
  been shown in Lemma \ref{lem:uvpositive}.\\
  Now for the verification of (\ref{9.3}) we let $0\le \varphi\in C_0^\infty(\bom\times [0,\infty))$
  be such that $\frac{\partial\varphi}{\partial\nu}=0$ on $\pO\times (0,\infty)$ and fix $T>0$ such that
  $\varphi\equiv 0$ in $\Omega\times [T,\infty)$.
  Then an application of Lemma \ref{lem16} shows that
  \begin{align}\label{eq:uepsforlimit}
 	c_1 \int_0^T \io &
	\veps^q |\nabla \ueps^\frac{p}{2}|^2 \varphi 
	+ c_2 \int_0^T\io \bigg| \ueps^{\frac p2}\nabla \ve^{\frac q2} 
	- \frac{(1-p)χ+2q}{2(pχ+1-q)}\ve^{\frac q2}\nabla \ue^{\frac p2}\bigg|^2φ\nn\\
	&=
	-\int_0^T \io \ueps^p \veps^q \varphi_t
	- \io u_0^p v_0^q \varphi(\cdot,0) 
	  +\frac{2pχ}q \int_0^T \io \ueps^\frac{p}{2} \veps^q \nabla\ueps^\frac{p}{2} \cdot\nabla\varphi \nn\\
	&  - \Big(1-\frac{p\chi}{q} \Big) \int_0^T \io \ueps^p \veps^q \Delta\varphi 
	  + q \int_0^T \io \ueps^p \veps^q \varphi 
	- q \int_0^T \io \frac{\ueps^{p+1}\veps^{q-1}}{1+ε\ueps}  \varphi,
  \end{align}
  where $c_1:=\frac{4(1-p)q-4q^2-p(1-p)^2χ^2}{pq(pχ+1-q)}$ and $c_2:=\frac{4(pχ+1-q)}{q}$ are positive, and 
  where employing \eqref{eq:upvqL1} we see that
  \be{77.9}
	- \int_0^T \io \ueps^p \veps^q \varphi_t
	\to - \int_0^T \io u^p v^q \varphi_t
  \ee
  and
  \be{77.10}
	q \int_0^T \io \ueps^p \veps^q \varphi
	\to q\int_0^T \io u^p v^q \varphi_t
  \ee
  as well as
  \be{77.11}
	- \Big(1-\frac{p\chi}{q}\Big) \int_0^T \io  \ueps^p \veps^q \Delta\varphi
	\to -\Big(1-\frac{p\chi}{q}\Big) \int_0^T \io u^p v^q \Delta\varphi
  \ee
  as $\eps=\eps_j\searrow 0$. Moreover, 
  $\ue^{\f p2}\ve^q\na\ue^{\f p2}= (\ue^{\f p2}\ve^{\f q2})(\ve^{\f q2}\na\ue^{\f p2})$, so that \eqref{eq:up2vq2L2} 
  and \eqref{eq:vq2naup2weakL2} ensure that
  \be{77.12}
 	\frac{2pχ}q \int_0^T \io 
	\ueps^\frac{p}{2} \veps^q \nabla\ueps^\frac{p}{2} \cdot\nabla\varphi \to \frac{2pχ}q \int_0^T \io 
	u^\frac{p}{2} v^q \nabla u^\frac{p}{2} \cdot\nabla\varphi
  \ee
  as $ε=ε_j\searrow0$. Combining $\eqref{eq:upp1vqm1L1}$ and Lebesgue's dominated convergence theorem, we furthermore 
  obtain
  \be{77.13}
 	- q \int_0^T \io \frac{\ueps^{p+1}\veps^{q-1}}{1+ε\ueps}  \varphi \to - q \int_0^T \io u^{p+1}v^{q-1}  \varphi
  \ee
  as $ε=ε_j\searrow0$. In view of a standard argument based on lower semicontinuity
  of the norm in $L^2(\Omega\times (0,T))$ with respect to weak convergence it follows from \eqref{eq:vq2naup2weakL2} 
  and \eqref{eq:up2navq2weakL2} and nonnegativity of $c_1$ and $c_2$ that inserting \eqref{77.9}, \eqref{77.10}, 
  \eqref{77.11}, \eqref{77.12} and \eqref{77.13} into \eqref{eq:uepsforlimit} yields 
  \begin{align*}
    	c_1 \int_0^T \io v ^q |\nabla u^\frac{p}{2}|^2 \varphi
   	&+ c_2 \int_0^T\io \bigg| u ^{\frac p2}\nabla v^{\frac q2} 
	- \frac{(1-p)χ+2q}{2(pχ+1-q)}v^{\frac q2}\nabla u^{\frac p2}\bigg|^2φ\nn\\
	&\le 	-\int_0^T \io u^p v^q \varphi_t - \io u_0^p v_0^q \varphi(\cdot,0) 
	+\frac{2pχ}q \int_0^T \io u^{\frac{p}{2}} v^q \nabla u ^{\frac{p}{2}} \cdot\nabla\varphi \nn\\
	 & - \Big(1-\frac{p\chi}{q} \Big) \int_0^T \io u ^p v ^q \Delta\varphi + q \int_0^T \io u ^p v ^q \varphi 
	- q \int_0^T \io u ^{p+1}v ^{q-1}  \varphi,
  \end{align*}
  which is equivalent to (\ref{9.3}) and thus completes the proof.
\qed
Our main result thereby becomes evident.\abs
\proofc of Theorem \ref{theo14}.\quad
  We only need to combine Lemma \ref{lem7} with Lemma \ref{lem77}.
\qed
%
%
%
%
%
%
\mysection{Appendix: A Poincar\'e inequality in non-convex domains}
Our proof of Lemma \ref{lem:iolnue} relies on Lemma \ref{lem:poincareforlog}, which in its original formulation in \cite[Lemma 4.2]{taowin_persistenceofmass} requires convexity of the domain due to the version of Poincar\'e's inequality (\cite[Corollary 9.1.4]{jost_pde}) used. 
In this appendix we state 
this Poincar\'e inequality without any such convexity condition, and since we could not find any 
reference to this in the literature, we briefly outline an argument.
Here and in the following, by $u_X$ we denote the average $\f1{|X|}\io u(x)dx$ for $u\in L^1(\Om)$ and any measurable set $X\subset \Om$ with positive measure.
\begin{lem}\label{lem:poincare}
  Let $\Om\sub ℝ^n$ be a bounded with smooth boundary, and let $δ>0$ and $p\in[1,\infty)$. 
  Then there exists $C=C(\Om,δ,p)$ with the property that for all $u\in W^{1,p}(\Om)$,
  \[
 	\kl{\io |u-u_B|^p}^{\f1p} \leq C(\Om,δ,p) \kl{\io |Du|^p}^{\f1p}
  \]
  holds for any measurable set $B\subset \Om$ with $|B|=δ$. 
\end{lem}
A derivation of this can be based on the following.
\begin{lem}\label{lem:replacejostlem913}
  Let $\Om\sub ℝ^n$ be a bounded domain with smooth boundary, and let $δ>0$. 
  Then there exists $C>0$ such that for every measurable set $B\subset \Om$ with $|B|=δ$, 
  and for each $u\in W^{1,1}(\Om)$ we have 
  \begin{equation}\label{eq:replacejostlem913}
   	|u(x)-u_B| \leq C\io \f{|Du(z)|}{|x-z|^{n-1}} dz
  \end{equation}
  for almost every $x\in \Om$.
\end{lem}
\proof
We follow the proof of \cite[Theorem 10]{hajlasz}, where \eqref{eq:replacejostlem913} is shown for $B=\Om$, and indicate necessary changes. With $B_0$ being a certain ball in $\Om$, defined as in the proof of \cite[Theorem 10]{hajlasz}, in \cite[(14)]{hajlasz} it is shown that there is $c_1>0$ such that for every $u\in W^{1,1}(\Om)$ and almost every $x\in \Om$ 
\begin{equation}\label{eq:hajlasz14}
 |u(x)-u_{B_0}|\leq c_1\io \f{|\na u(z)|}{|x-z|^{n-1}} dz. 
\end{equation}
Whereas the first summand in the right-hand side of 
\begin{equation}\label{eq:dreiecksungl}
 |u(x)-u_B|\leq |u(x)-u_{B_0}| + |u_{B_0}-u_B| 
\end{equation}
is immediately covered by \eqref{eq:hajlasz14}, as to the second we observe that, again by \eqref{eq:hajlasz14}, 
\begin{align}\label{eq:uB0minusuB}
 |u_{B_0}-u_B| & \leq \f1{|B|}\int_B |u_{B_0}-u(y)| dy\nn \\
 &\leq \f {c_1}{|B|} \int_B\io \f{|\na u(z)|}{|y-z|^{n-1}} dz dy\nn\\ 
 &\leq \f{c_1}{|B|} \io |\na u(z)| \int_B\f{1}{|y-z|^{n-1}}dydz.
\end{align}
Here we use that with some $c_2=c(n)$, for all $z\in ℝ^n$ and any measurable $E\subsetℝ^n$, $\int_E \f{dy}{|y-z|^{n-1}} \leq c_2 |E|^{\f1n}$ (\cite[(13)]{hajlasz}) and that $1\leq \f{(\diam\Om)^{n-1}}{|x-z|^{n-1}}$ for any $x,z\in \Om$. With these observations, \eqref{eq:uB0minusuB} turns into 
\begin{equation}\label{eq:poincareprooffinal}
 |u_{B_0}-u_B| \leq \f{c_1c_2}{|B|}|B|^{\f1n} \io |\na u(z)| dz \leq c_1c_2|B|^{\f1n-1} (\diam \Om)^{n-1} \io \f{|\na u(z)|}{|x-z|^{n-1}} dz. 
\end{equation}
Noting that $|B|^{\f1n-1}=δ^{\f1n-1}$ and combining \eqref{eq:hajlasz14} and \eqref{eq:poincareprooffinal} with \eqref{eq:dreiecksungl} proves \eqref{eq:replacejostlem913}.
\qed
\proofc of Lemma \ref{lem:poincare}. \quad
  For convex domains, this is exactly Corollary 9.1.4 of \cite{jost_pde}, 
  which follows from \cite[Lemma 9.1.3]{jost_pde} and \cite[Lemma 9.1.2]{jost_pde}, the latter of which 
  (a continuity property of the Riesz potential operator) poses no convexity condition on $\Om$. 
  As replacement of the former, in the case of general $\Omega$ we now rather rely on Lemma \ref{lem:replacejostlem913}.
\qed
\begin{rem}
  In Lemma \ref{lem:replacejostlem913} (and hence in Lemma \ref{lem:poincare}), for the domain it is actually sufficient 
  to be (bounded and) a John domain, instead of having smooth boundary. 
  In particular, any bounded domain satisfying the interior cone condition is admissible in these lemmata. 
  For details concerning this, we once more refer the reader to \cite{hajlasz}. 
\end{rem}
\begin{rem}
  With Lemma \ref{lem:poincare}, it is also possible to remove the convexity condition on the domain 
  in \cite{taowin_persistenceofmass}. 
\end{rem}
%
%
%
%
%
%
%
%
%
%
%
%
{\bf Acknowledgement.} \quad
The authors acknowledge support of the DFG within the project ``Analysis of chemotactic cross-diffusion in complex frameworks''. 

\begin{thebibliography}{1}
%
\bibitem{BBTW}
  \sc Bellomo, N., Bellouquid, A., Tao, Y., Winkler, M.:
  \it Toward a mathematical theory of Keller-Segel models of pattern formation in biological tissues. 
  \rm Math.~Models Methods Appl.~Sci. {\bf 25}, 1663-1763 (2015)
\bibitem{biler}
  \sc Biler, P.: \it Global solutions to some parabolic-elliptic systems of chemotaxis. 
  \rm Adv.~Math.~Sci.~Appl. {\bf 9} (1), 347-359 (1999)
\bibitem{diperna_lions}
  \sc Di Perna, R.-J., Lions, P.-L.: \it
  On the Cauchy problem for Boltzmann equations: Global existence and weak stability.
  \rm Ann. Math. {\bf 130}, 321-366 (1989)
\bibitem{fujie_bounded}
  \sc Fujie, K.: \it Boundedness in a fully parabolic chemotaxis system with singular sensitivity.
  \rm J.~Math.~Anal.~Appl. {\bf 424} 675-684 (2015)
\bibitem{fujie_senba}
  \sc Fujie, K., Senba, T.: 
  \it Global existence and boundedness in a parabolic-elliptic Keller-Segel system with general sensitivity. 
  \rm Discr.~Cont.~Dyn.~Syst.~B {\bf 21}, 81-102 (2016)
\bibitem{fujie_senba_NON}
  \sc Fujie, K., Senba, T.:  \it Global existence and boundedness of radial solutions to a two dimensional 
  fully parabolic chemotaxis system with general sensitivity. 
  \rm Nonlinearity {\bf 29}, 2417–2450 (2016)
\bibitem{hajlasz}
  \sc Haj{\l}asz, P.: \it Sobolev inequalities, truncation method, and {J}ohn domains.
  \rm  In: {\em Papers on analysis}, {\em Rep. Univ.
  Jyv\"askyl\"a Dep. Math. Stat.} {\bf 83}, 109-126 (2001)
\bibitem{herrero_velazquez}
  \sc Herrero, M.~A., Vel\'azquez, J.~J.~L.:
  \it A blow-up mechanism for a chemotaxis model.
  \rm Ann.~Scuola Normale Superiore Pisa Cl.~Sci. {\bf 24}, 633-683 (1997)
\bibitem{hillen_painter2009}
  \sc Hillen, T., Painter, K.J.: \it A User's Guide to PDE Models for Chemotaxis.
  \rm J.~Math.~Biol. {\bf 58} (1), 183-217 (2009)
\bibitem{jost_pde}
  \sc Jost, J.: \it Partial differential equations. 
  \rm Graduate Texts in Mathematics, Vol. 214. Springer, New York, second edition, 2007
\bibitem{lankeit}
  \sc Lankeit, J.: 
  \it A new approach toward boundedness in a two-dimensional parabolic chemotaxis system with singular sensitivity. 
  \rm Math.~Meth.~Appl.~Sci. {\bf 39}, 394-404 (2016)
\bibitem{luckhaus_sugiyama_ARMA}
  \sc Luckhaus, S., Sugiyama, Y., Vel\'azquez, J.J.L.: 
  \it Measure valued solutions of the 2D Keller-Segel system. 
  \rm Arch.~Rat.~Mech.~Anal. {\bf 206}, 31-80 (2012)
\bibitem{mizukami_yokota}
  \sc Mizukami, M., Yokota, T.:
  \it A unified method for boundedness in fully parabolic chemotaxis systems with signal-dependent sensitivity.
  \rm preprint, 2017. https://arxiv.org/abs/1701.02817
\bibitem{nagai_senba1998}
  \sc Nagai, T., Senba, T.: 
  \it Global existence and blow-up of radial solutions to a parabolic-elliptic system of chemotaxis.
  \rm Adv.~Math.~Sci.~Appl. {\bf 8}, 145-156 (1998)
\bibitem{rosen}
  \sc Rosen, G.: \it Steady-state distribution of bacteria chemotactic toward oxygen.
  \rm Bull.~Math.~Biol. {\bf 40}, 671-674 (1978)
\bibitem{simon}
  \sc Simon, J.: \it Compact sets in the space {$L^p(0,T;B)$}.
  \rm Ann. Mat. Pura Appl. {\bf 146}, 65-96 (1987)
\bibitem{stiwi}
  \sc Stinner, C., Winkler, M.: \it Global weak solutions in a chemotaxis system with large singular sensitivity.
  \rm Nonlinear Analysis: Real World Applications {\bf 12}, 3727-3740 (2011)
\bibitem{taowin_persistenceofmass}
  \sc Tao, Y., Winkler, M.: \it Persistence of mass in a chemotaxis system with logistic source.
  \rm J.~Differential Eq. {\bf 259}, 6142-6161 (2015)
\bibitem{tello_win_Pisa}
  \sc Tello, J.I., Winkler, M.: 
  \it Reduction of critical mass in a chemotaxis system by external application of a chemoattractant. 
  \rm Ann.~Sc.~Norm.~Sup.~Pisa Cl.~Sci. {\bf 12}, 833-862 (2013)
\bibitem{win_aggregationvs}
  \sc Winkler, M.: \it Aggregation vs. global diffusive behavior in the higher-dimensional Keller-Segel model.
  \rm J.~Differential Eq. {\bf 248}, 2889-2905 (2010)
\bibitem{win_MMAS}
  \sc Winkler, M.: \it Global solutions in a fully parabolic chemotaxis system with singular sensitivity.
  \rm Math.~Meth.~Appl.~Sci. {\bf 34}, 176-190 (2011)
\bibitem{win_JMPA}
  \sc Winkler, M: \it Finite-time blow-up in the higher-dimensional parabolic-parabolic Keller-Segel system.
  \rm J.~Math.~Pures Appl. {\bf 100}, 748-767 (2013), {\tt arXiv:1112.4156v1}
\bibitem{win_SIMA}
  \sc Winkler, M.: \it Large-data global generalized solutions in a chemotaxis system with tensor-valued sensitivities.
  \rm SIAM J.~Math.~Anal. {\bf 47}, 3092-3115 (2015)
\bibitem{zhao_zheng_sining_JMAA2016}
  \sc Zhao, X., Zheng, S.: 
  \it Global boundedness of solutions in a parabolic-parabolic chemotaxis system with singular sensitivity. 
  \rm J.~Math.~Anal.~Appl. {\bf 443}, 445-452 (2016)
%
\end{thebibliography}
\end{document}